\documentclass{article}
\usepackage{geometry}
\usepackage{bm,epsfig,ifthen}
\usepackage{amsfonts,amssymb}
\usepackage[english]{babel}
\usepackage[latin1]{inputenc}
\usepackage{graphicx}
\usepackage{amssymb,amsmath}
\usepackage{pstricks}

\newcommand{\field}[1]{\mathbb{#1}}

\title{Statistical behavior of adaptive multilevel splitting algorithms in simple models}

\author{Joran Rolland \footnote{INLN, UMR 7335 CNRS, UNSA, 1361 route des lucioles, 06560 Valbonne, France} \footnote{joran.rolland@inln.cnrs.fr} , Eric Simonnet$^*$ \footnote{eric.simonnet@inln.cnrs.fr}}

\date{Submitted to J. Comp. phys. : 12/05/2014, Accepted: 3/12/2014.}
\begin{document}

\maketitle

\begin{flushleft}Keywords :Multilevel splitting, Stochastic processes, Numerical methods,
 phase transition\end{flushleft}\begin{flushleft}PACS:  05.40.-a,  02.60.-x\end{flushleft}

\begin{abstract}
Adaptive multilevel splitting algorithms have been introduced rather recently for estimating tail distributions
in a fast and efficient way. In particular, they can be used for computing the so-called reactive trajectories
corresponding to direct transitions from one metastable state to another.
The algorithm is based on successive selection-mutation steps performed on the system in a controlled way.
It has two intrinsic parameters, the number of particles/trajectories and
the reaction coordinate used for discriminating good or bad trajectories.
We investigate first the convergence in law of the algorithm as a function of the
timestep for several simple stochastic models.
Second, we consider the average duration of reactive trajectories for which no theoretical predictions
exist.
The most important aspect of this work concerns some systems with two degrees
of freedom. They are studied in details as a function of the reaction coordinate in
the asymptotic regime where the number of trajectories goes to infinity.
We show that during phase transitions, the statistics of the algorithm deviate significatively from
known theoretical results when using non-optimal reaction coordinates.
In this case, the variance of the algorithm is peaking at the transition and the convergence of
the algorithm can be much slower than the usual expected central limit behavior.
The duration of trajectories is affected as well. Moreover, reactive trajectories do
not correspond to the most probable ones. Such behavior disappears when using the optimal reaction coordinate
called committor as predicted by the theory.
We finally investigate a three-state Markov chain which reproduces this phenomenon and show
logarithmic convergence of the trajectory durations.

\end{abstract}
\maketitle

\section{Introduction}

 Computing reactive trajectories between two metastable states $\mathcal{A}$ and $\mathcal{B}$ and
the associated crossing probability $\alpha$ is a central problem
to kinetic chemistry \cite{ha,jpcm,ffsf,jcp1,jcp2}  (Fig.~\ref{sketches} (a)). Such rare excursions arise in many
other fields, like climate science \cite{Ey}, geophysics \cite{dyn}, fluid dynamics \cite{BS,physrep}, probabilities and statistics \cite{cglp,poi}. However, the estimation of these transitions requires special care because it often
concerns extremely rare events with probabilities of order $O(10^{-10})$ and less.
For this reason, direct numerical simulations or Monte Carlo Markov Chains are out of the question.

In order to tackle this question, different approaches can be considered. Transition rates calculations go back to the Eyring--Kramers theory of mean first passage times \cite{ha} (Fig.~\ref{sketches} (a)).
In mathematics, Freidlin--Wentzell theory of large deviations allows one to compute instantons
corresponding to the most probable trajectories between one state to another \cite{ht}.
These trajectories are obtained by minimizing some action of the model and more importantly
they are achieved in the so-called zero-noise limit of the system.
One of the problems is that for large-dimensional systems, the minimization problem often
becomes ill-posed and numerically very difficult to solve.
Moreover, large deviations often require particular conditions of applications
and in fact might not even exist. It is known for instance that Freidlin--Wentzell theory
cannot be applied when local attractors are non-isolated (see \cite{bt_niso}),
a situation which happens in 2-D turbulence in the so-called Eulerian limit \cite{physrep}. Other,
more general, approaches are needed and Adaptive Multilevel Splitting (AMS) algorithms
\cite{cglp,poi,cg07,pdm} appear to be promising.

AMS is a method which falls in the class of multilevel splitting algorithm
and more generally Forward Flux Sampling (FFS) \cite{jpcm,ffsf}.
Let $X$ be a random variable taking values in some space ${\cal E}$, and some application
$$
Q: {\cal E} \to \field{R}.
$$
Multilevel splitting algorithm aims at estimating
\begin{equation}\label{a0}
\alpha \equiv \field{P}(Q(X) > l),
\end{equation}
by decomposing $\alpha$ as
\begin{equation}\label{product}
\alpha = \prod_k \field{P}(Q(X) > l_{k+1}|Q(X) > l_k),~-\infty < \cdots < l_k < l_{k+1} < \cdots +\infty.
\end{equation}
The idea is that the intermediate conditional probabilities are larger than $\alpha$
and thus should be estimated more easily.
To do so, one uses a splitting approach where the system is duplicated in $N$ parallel copies.
At each step, only the proportion of particles being above a given level is kept whereas the
others particles are simply killed. New particles are then created to keep $N$ constant during
the course of the algorithm.
The advantage of these approaches is the possibility to go
beyond the Freidlin--Wentzell theory. Firstly, they give access to the whole distribution of reactive trajectories. By contrast, the Freidlin--Wentzell
theory only provides the most probable trajectory. Secondly They can be
applied for very general systems and for finite-noise amplitude without
any theoretical restrictions apart from dealing with (strong) Markovian systems.

Classical multilevel algorithms correspond to the situation where the levels $l_k$ are
fixed a priori. It appears that most often, the choice of these levels is non-trivial
and requires a good knowledge of the system. For instance, it often happens that none
of the particles are able to reach a given level, leading to a premature ending of the algorithm.
The idea of adaptive multilevel
splitting algorithms is to choose these levels in an adaptive way during the course of
the algorithm.  The strategy is to retain at each step a fixed proportion of ``good" particles,
say $N-n$, and kill the $n$ remaining ``bad" particles, i.e. those having the lowest values of $Q$.
These $n$ particles are then replaced by a set of new particles which inherit some of the properties
of the successful $N-n$ particles: this is why this procedure is sometimes referred
to as a {\it selection-mutation} step.

In practice, such step can follow various strategies. A natural and simple
one is to restart the $n$ trajectories from different initial conditions given by the
$N-n$ good particles: a detailed description of how to do that is given in the first section.
The important point is that the different $l_k$ are now random which complicates theoretical
analysis quite a bit.

The situation where one is considering a strong Markov process $(X_t)_{t \in \field{R}}$, that is
the variable $X$ above depends on some parameter in a nontrivial way, is called in the
literature the {\it dynamical case}. In other words, one is working at the level of trajectories
(still named particles). This is the situation we will consider in this work.
On the other hand, when the random variable $X \in {\cal E}$ does not depend on time,
the problem is referred to as the {\it static case} \cite{pdm}.

In this paper, we focus on the more specific problem of finding an estimate of the crossing probability $\alpha$
in the dynamical case. We want to find the probability $\alpha$
that a trajectory starting from some initial condition $x$ reaches a set ${\cal B}
\in \field{R}^d$ before some set ${\cal A} \in \field{R}^d$ with ${\cal A} \cap {\cal B} = \emptyset$
(see Fig.~\ref{sketches} (a)). It is defined as
\begin{equation}\label{a1}
\alpha \equiv \field{P}(\tau_{\cal B}(x) < \tau_{\cal A}(x)),~{\rm with}~\tau_{{\cal A},{\cal B}} =
\inf_t \{ X_t \in {\cal A},{\cal B}; X_0 = x \}.
\end{equation}
The probability $\alpha$ seen as a function of $x$ is called the {\it committor} or {\it equilibrium potential}
in mathematics (see \cite{jpcm,VE}). To make the connection with multilevel splitting algorithm and
the quantity (\ref{a0}), one defines first
the so-called {\it reaction coordinate}:
$$
\Phi: \field{R}^d \to \field{R}
$$
and the application $Q$
from the space of trajectories ${\cal E}$ which start from $x$ to $\field{R}$ is defined as
$$
Q:{\cal E} \to \field{R}, (X_t)_{t \in [0,\tau_{\cal A}]} \mapsto \sup_{t \in [0,\tau_{\cal A}]}
\Phi(X_t), ~X_0 = x.
$$
The reaction coordinate therefore measures how good the trajectories are by telling how far
they can escape from the set ${\cal A}$.
Finally, when ${\cal B}$ is such that ${\cal B} = \Phi^{-1}(b)$, one has:
$$
\field{P}(\tau_{\cal B}(x) < \tau_{\cal A}(x)) = \field{P}(\sup_{t \in [0,\tau_{\cal A}]} \Phi(X_t) > b),~
X_0 = x.
$$

Other adaptive algorithms have been proposed \cite{jcp1}. Some include strategies to limit the simulation time
when trajectories are trapped in a local metastable state \cite{jcp2}, but one loses
the ability to compute the crossing probability.
AMS distinguishes itself in two aspects: first it allows one to compute all the properties of
the ``reaction'' \cite{cglp};
second, the convergence of the crossing probability estimate as the number of copies $N$
goes to infinity has been studied mathematically \cite{poi,cg07,pdm}.
In fact, most of the known theoretical results assume that the intermediate probabilities in (\ref{product})
are sampled exactly: in this case, one speaks of {\it perfect} or {\it ideal} AMS.

In that case, one can predict the behavior
of the variance and the bias on the estimation of $\alpha$ \cite{cg07,pdm} which turn out to
scale like $1/\sqrt{N}$ for the variance and $1/N$ for the bias.
More precisely, in the static perfect case and when one eliminates only $n=1$ particles
at each step (the so-called last particle method),
the total number of iterations $K$ when the algorithm stops has been shown \cite{poi} to follow:
\begin{equation}\label{poisson}
K \sim {\rm Poisson}(-N \log \alpha)
\end{equation}
(see Fig~\ref{dist_num} (a)).
In the 1-D dynamical case, assuming the trajectories are continuous,
a central limit theorem has been established by \cite{cg07}, taking the form:
\begin{equation}\label{clt}
\begin{array}{l}
\sqrt{N} \left(\alpha - \hat \alpha\right) \xrightarrow[N \to \infty]{{\cal D}}
{\cal N}(0,\sigma^2),~ \hat \alpha \xrightarrow[N \to \infty]{a.s.} \alpha
\end{array}
\end{equation}
(see Fig.~\ref{dist_num} (b)), where
$$
\hat \alpha = (1-\rho) (1-q)^{K}, ~\sigma^2 = \alpha^2 \left(K_0 \frac{q}{1-q} + \frac{\rho_0}{1-\rho_0} \right),
$$
with $(\rho_0,K_0)$ such that $\alpha = (1-\rho_0)(1-q)^{K_0}$. Here, The quantile $q$ is assumed to
be independent on $N$ and corresponds to the prescribed
proportion of particles which are rejected at each step. For finite $N$, the quantile would be defined as
$q=n/N$. Moreover,
this result does not take into account discretized Markov (jump) processes.
Note that the two results (\ref{poisson}) and (\ref{clt}) are consistent with each other although at a formal level only (taking $q=1/N$).
In that case, the variance becomes in the limit of large $N$
\begin{equation}\label{variance1}
\sigma^2 = -\alpha^2 \log \alpha.
\end{equation}
More results can be obtained on the bias, for instance in the static perfect case
one can show that $\langle \hat \alpha \rangle - \alpha > 0$ which
decreases like $N^{-1}$ \cite{pdm}.
These results have been numerically confirmed  in the 1-D case
by \cite{cg07} and \cite{cglp} and to some extend for 2-D systems as well when the
system temperature is not too small.

The case where the intermediate probabilities (\ref{product}) are not sampled exactly yields to
serious theoretical difficulties and the statistical behavior of the algorithm in the dynamical
case $d > 1$ is unknown due to the complex dependency of the algorithm on the chosen
reaction coordinate $\Phi$. A good choice of $\Phi$ in fact yields near optimal variance, scaling
like $1/\sqrt{N}$, and nice convergence bias properties as well. In fact, using the committor function for
$\Phi$ can be shown to give optimal results \cite{STCO}.

For complex systems however, choosing a good reaction coordinate is a difficult issue.
Using the committor is out of question since it is precisely the quantity one would like to estimate.
The only exception is for very low dimensional system (typically less than $d \leq 3$),
where the committor can be computed
directly by other means, in which case AMS might lose part of its attractiveness.
Another open issue is the average duration of reactive trajectories $\tau$ found by the algorithm
for which no theoretical results exist.

The key results of this paper is that, even for systems with only two degree of freedom, one
can observe a bias of different amplitude and sign than what is predicted. More importantly,
the variance may be much larger than expected in the asymptotic
regime $N \to \infty$. In fact, in the last particle algorithm situation $n=1$, the number of iterations
can strongly deviate from the pure Poisson law (\ref{poisson}).
These phenomena become more apparent as the temperature goes to zero (large bias)
or during phase transitions (large variance).
It appears therefore that the choice of the reaction coordinate becomes highly critical \cite{jpcm,ffsf,cglp}.

The plan of this work is as follow. After describing in details the algorithm and
the models we consider in Section \S~\ref{S1}, we illustrate some convergence results for these models
in Section \S~\ref{S2}.  In Section \S~\ref{reac}, we discuss in details the quality criterion
for reaction coordinates. Finally, we introduce a simple model
for understanding the effect of a poorly-chosen reaction coordinate. We then quantify the rate of
convergence of the algorithm in that case.
The results are eventually discussed altogether in the conclusion (\S~\ref{conc}).

\section{AMS description and models \label{S1}}
In this section we provide a detailed description of the algorithm
(see also \cite{cglp,cg07,pdm}). In order to treat a slightly
more general and practical situation, we will allow the algorithm to start not from a single initial condition
but from a set of initial conditions distributed on an hypersurface ${\cal C}$ surrounding ${\cal A}$ (Fig.~\ref{sketches}(a)).
Let $\rho_{{\cal C}}$ be the restriction of the equilibrium probability measure of the process $(X_t)$ to
the hypersurface ${\cal C}$
so that the initial starting conditions are distributed according to $\rho_{{\cal C}}$.
The probability $\alpha$ now corresponds to the quantity
$$
\alpha \equiv \field{P}(\tau_{\cal B}(x) < \tau_{\cal A}(x)|x \sim \rho_{\cal C})
= \int_{\cal C} \field{P}(\tau_{\cal B}(x) < \tau_{\cal A}(x)) \rho_{\cal C}(x)\,{\rm d}x.
$$
We moreover consider very general sets ${\cal A}$ and ${\cal B}$ under the natural constraint that
${\cal A} \cap {\cal B} = \emptyset$. Without loss of generality, we will
assume that the reaction coordinate $\Phi$ satisfies
$$
\Phi: \field{R}^d \to \field{R},~ 0 \leq \Phi \leq 1,~\Phi(\partial {\cal A}) = 0,~\Phi(\partial {\cal B}) = 1.
$$
As explained in the introduction, this application indeed monitors the position in phase space relative to
the sets ${\cal A}$ and ${\cal B}$.

\subsection{The algorithm \label{ss3}}

\begin{figure}
\includegraphics[width=8cm]{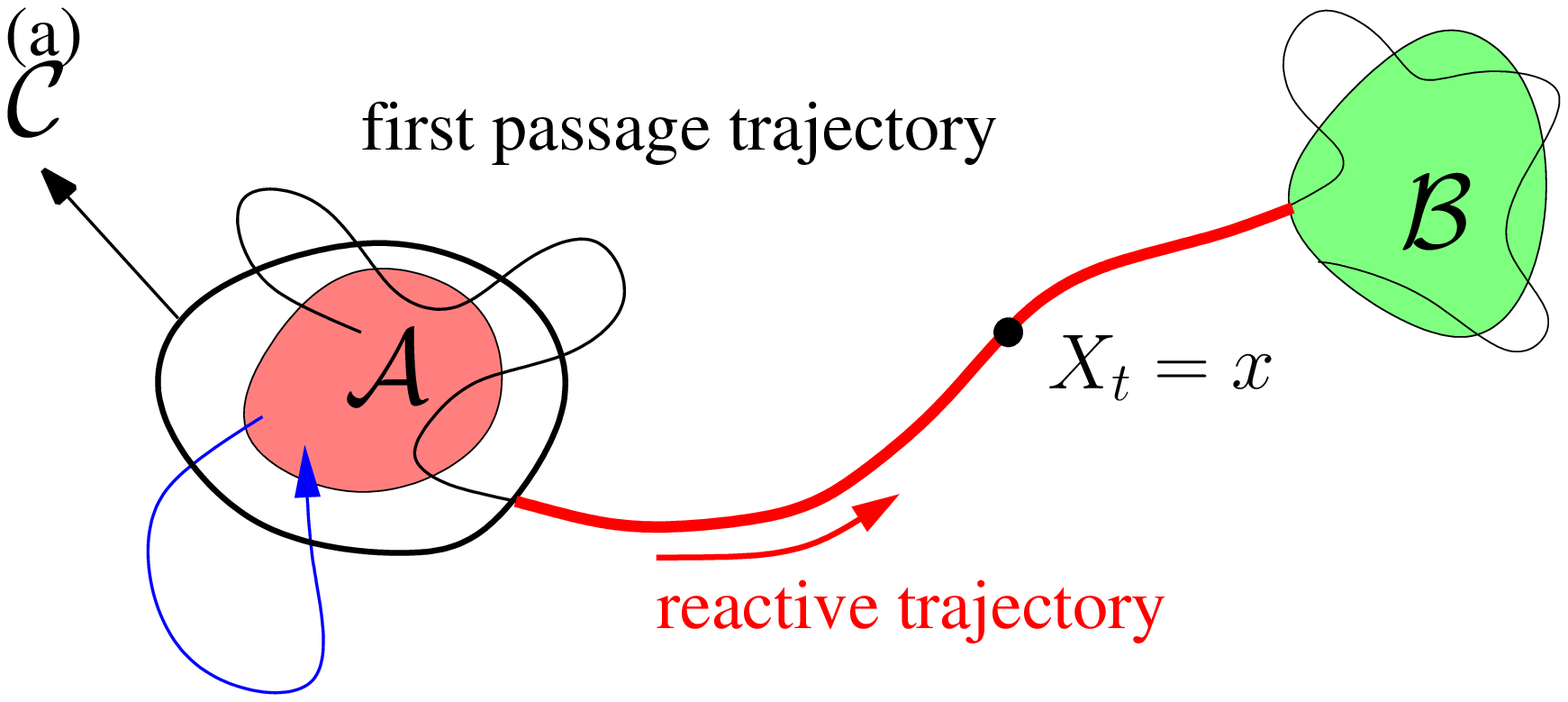}
\includegraphics[width=6.5cm]{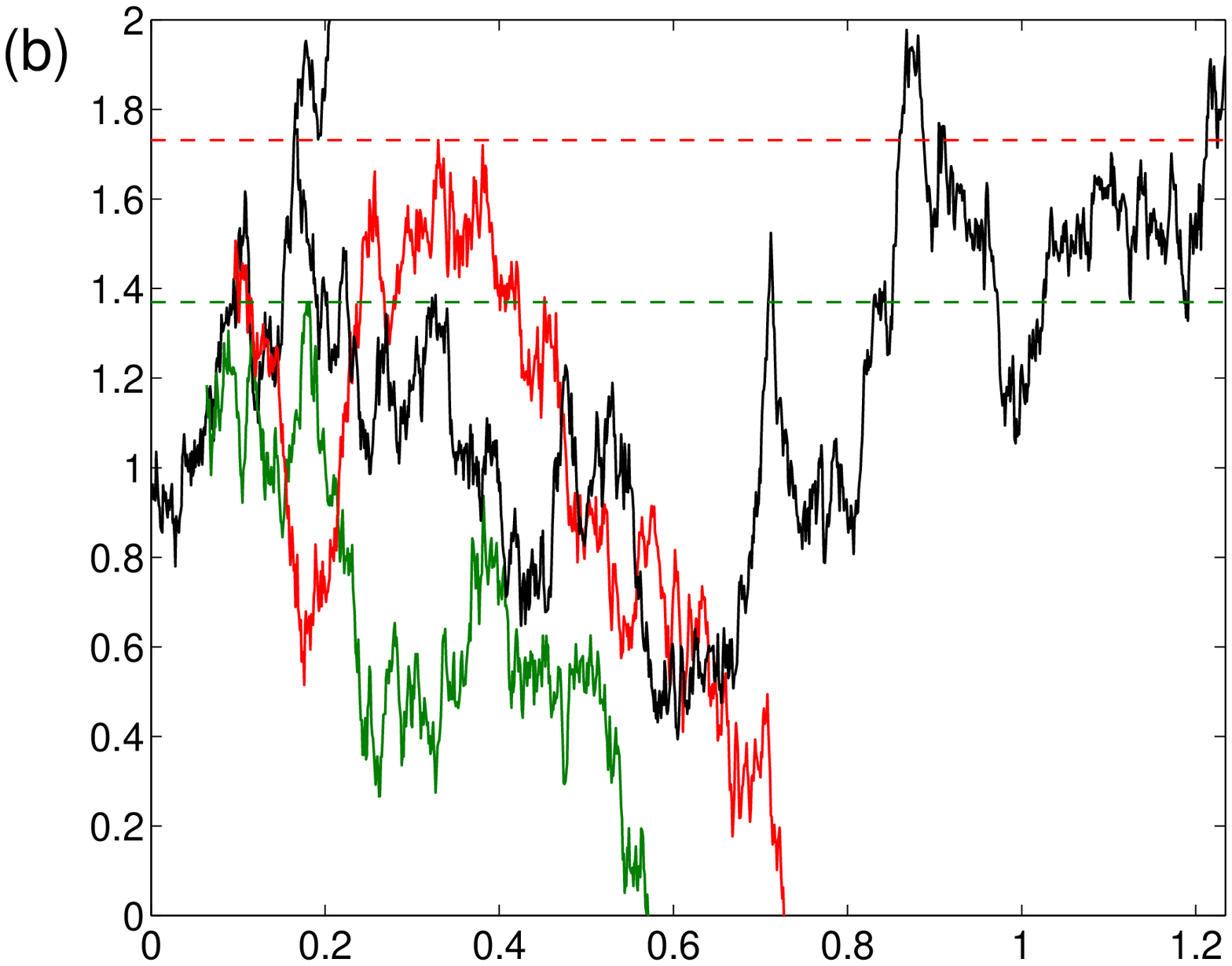}
\caption{(a): a first passage trajectory from $\mathcal{A}$ to $\mathcal{B}$ together with the hypersurface
$\mathcal{C}$,
is shown in black and red with its portion in red corresponding
to a reactive trajectory. A non-reactive trajectory is shown in blue.
(b): example of AMS on the Brownian drift (See \S~\ref{brdr}), $\mu=0.3$, $N=2$, $\mathcal{A}=0$, $\mathcal{B}=2$.
The initial condition is $1$.
Black lines correspond to reactive trajectories, coloured lines to suppressed non reactive trajectories, dashed lines correspond to the adapted levels. Some non-reactive trajectories going directly to $0$ are not displayed for readability.}
\label{sketches}
\end{figure}
The number of iterations of the algorithm is denoted by $k$.
We give a concrete illustration of these steps for a Brownian drift with $N=2$. (see Eq.~(\ref{com_drift})). The algorithm can
be decomposed into the following steps:
\begin{enumerate}
\item{\it Initial step $k=1$}\\
Fix the total number of particles $N$, the quantile $n$ of particles killed at each step verifying $1\leq n \leq N-1$ and the reaction coordinate $\Phi$.
Define an hypersurface ${\cal C}$.
Generate $N$ independent identically distributed (i.i.d.) trajectories $X_t^i$, $i = 1,\cdots,N$ starting from $N$
i.i.d. initial conditions distributed on $\mathcal{C}$ according to $\rho_{\mathcal{C}}$.
Wait until all of them either reach $\mathcal{A}$ or $\mathcal{B}$
(see Fig. \ref{sketches} (b): black curve).

\item{\it Selection-mutation}
\begin{itemize}
\item
Compute the $N$ values
$$
Q_i \equiv \sup_{t \in [0,\tau_{\cal A}^i]} \Phi(X_t^i), ~i = 1,\cdots,N.
$$
Sort the $N$ values so that $Q_{(1)} \leq Q_{(2)} \cdots \leq Q_{(N)}$ and set
$$
l_k \equiv Q_{(n)}.
$$
Kill the $n$ trajectories corresponding to $Q_{(1)},\cdots Q_{(n)}$
(see Fig. \ref{sketches} (b): green and red curves), with $Q_{(m)}$ the $m$th ordered value of the $\{Q_i\}$.
\item In order to keep the total number $N$ of trajectories fixed, one needs to restart the $n$ trajectories
which have been killed. They are generated by branching them on the $N-n$ good trajectories.
This branching (mutation) procedure is the following: for each trajectory
being killed, with index $i$ say and maximum level $Q_i$, we associate a trajectory chosen randomly with uniform probability $1/(N-n)$
among the $N-n$ good trajectories, call its index $i^\star$.
Compute a new initial condition $x^\star$ corresponding to the intersection with
the isosurface $\Phi^{-1}(Q_i)$ and
the good trajectory with index $i^\star$. This defines a time $t^\star$ such that $x^\star = X^{i^\star}_{t^\star}$.
In practice, since the system is discretised, the trajectories are not continuous, the best approach is to
take instead
$$
t^\star \equiv \inf_t \{\Phi(X^{i^\star}_t) \geq Q_i \},~x^\star \equiv X^{i^\star}_{t^\star}.
$$
Do this for the $n$ bad trajectories yielding $n$ initial conditions $x^\star_i$.
\item Restart the $n$ trajectories from the $n$ new initial conditions $x_i^\star$
until they all reach either ${\cal A}$
or ${\cal B}$ (Fig.~\ref{sketches} (b), green, red and second black curves).
$$
k \to k + 1\,.
$$

\end{itemize}
\item  The algorithm stops after a (random) number of iterations say $k$
when $N-n$ trajectories have reached ${\cal B}$ before ${\cal A}$. In such a
case they are in general a proportion $r \geq N-n$ of trajectories having done so.
Note that in the particular case of choosing $n=1$ (the so-called last particle method) then $r=N$.
Call the total number of iterations $K \equiv k-1$.
\end{enumerate}

The numbers $r$ and $K$ are random with different values for different realisation of the
algorithm. They yield an estimate of the crossing probability \cite{cglp,cg07} :
\begin{equation}
\label{alpha_algo}
\hat{\alpha}= \left(\frac{r}{N}\right) \left(1-\frac{n}{N}\right)^{K}\,.
\end{equation}
By a simplification of language both the estimate $\hat{\alpha}$ and $\alpha$ will be termed crossing probability.

\begin{figure}
\centerline{\includegraphics[width=6cm]{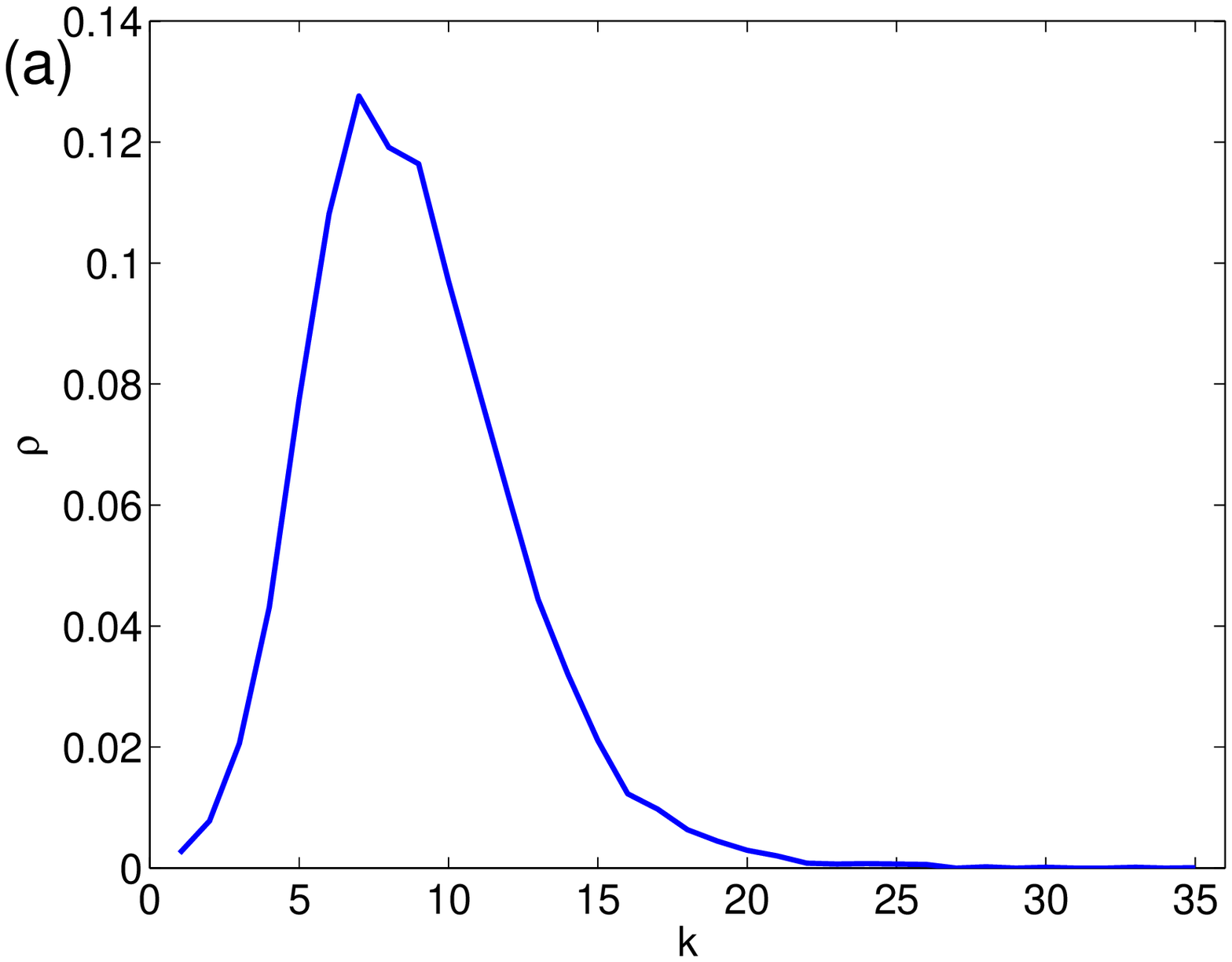}
\includegraphics[width=6cm]{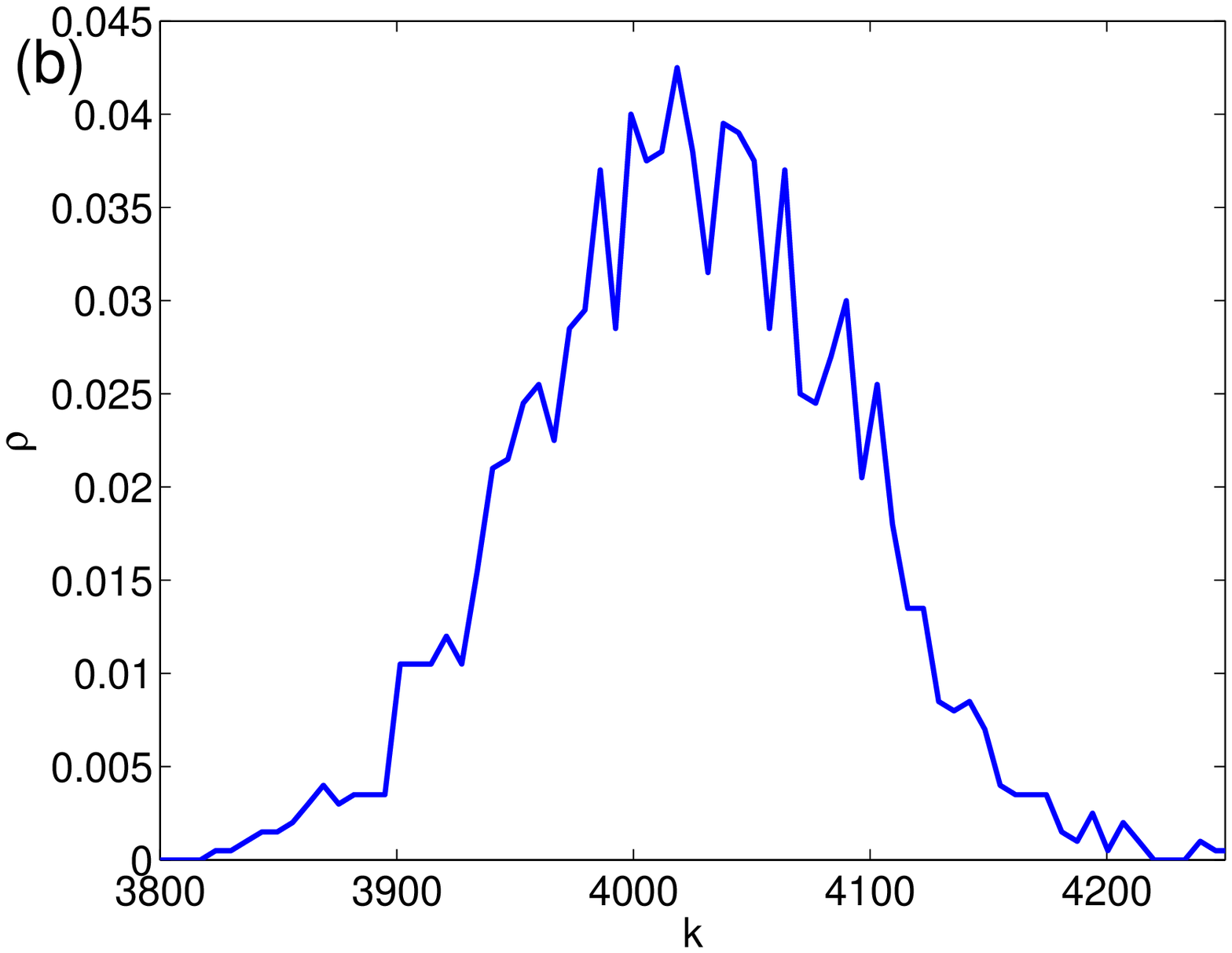}
\includegraphics[width=6cm]{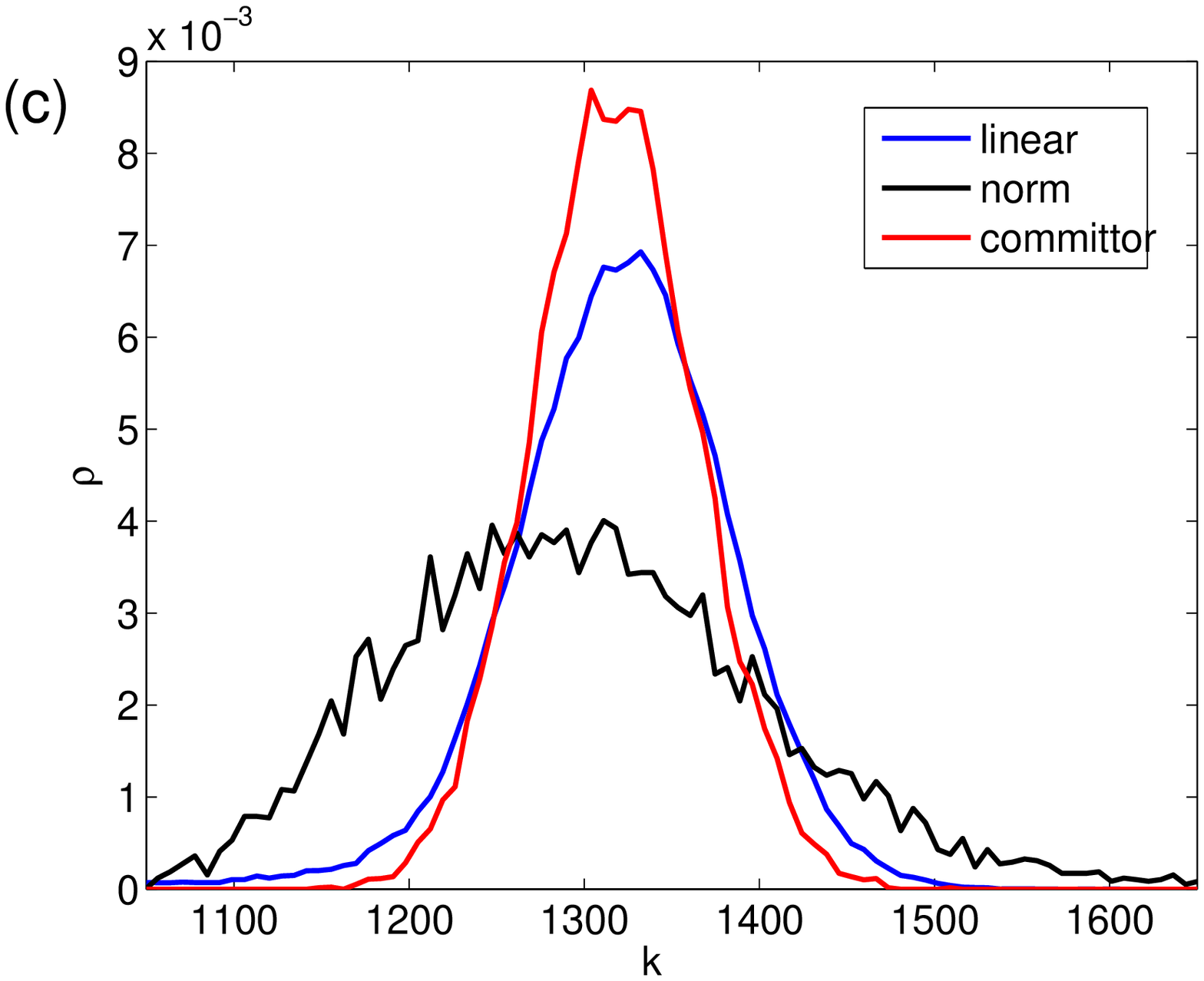}}
\caption{Examples of distributions of the number of iterations of the algorithm applied
to the triple-well model for several regimes of parameters. (a) :
$N=3$ clones and $\beta=0.05$ (histogram). (b) $N=1000$ clones, $\beta=1$ (histogram), (c) : $N=100$ clones,
$\beta=5$, two reaction coordinates (Eq.~(\ref{coord_2c})) and reaction
coordinate based on the committor (Fig.~\ref{comm2d}) (pdf).}
\label{dist_num}
\end{figure}

\subsection{Models}
In this part, we present the models used for investigating the detailed statistical behavior of the algorithm.
These models are all of the form of a stochastic differential equation \cite{gar,VK}:
\begin{equation}
\frac{dX}{dt}=F(X)+\sqrt{\frac{2}{\beta}}\eta\,,\,
X \in \mathbb{R}^n\,,\, \langle \eta_j(t)\eta_{j'}(t')\rangle=\delta(t-t')\delta_{jj'}\,,
\label{eqpr}
\end{equation}
with $\beta=1/T$, the inverse of the temperature.
This equation is referred to as the over-damped Langevin or gradient systems when the force $F$ derives from a potential $V$, i.e.
$$
F = -\nabla_X V.
$$
In the following, we will consider such systems. The advantage is that their dynamics are
directly given by the potential $V$.
To fix the idea, we assume that these systems exhibit some form of metastability.
The sets $\mathcal{A}$ and $\mathcal{B}$ correspond to some neighbourhood of two stable fixed points of $F$,
in which $X$ spend most of the time. The hypersurface $\mathcal{C}$
closely surrounds $\mathcal{A}$ as shown in Fig. \ref{sketches} (a).
Of interest is the committor function (see Eq.~(\ref{a1})).
For diffusive processes of gradient form, it can be shown to solve the backward Fokker-Planck equation:
\begin{equation}\label{comm}
\mathbf{F}.\mathbf{\nabla}q+\frac{1}{\beta}\Delta q = 0
\,,\, \forall x\in \partial \mathcal{A}\,,\, q(x)=0  \,,\,  \forall x\in \partial \mathcal{B}\,,\,
q(x)=1\,.
\end{equation}
Beside giving the probability of crossing, the committor has several interesting properties.
Indeed, one can extract the flux of reactive trajectory from $q$ (see \cite{VE}),
which indicates the most likely paths that the system can follow between $\mathcal{A}$ and $\mathcal{B}$.
More importantly, as it has been noted \cite{jpcm,cglp} and as we will show,
the function $\Phi_c \equiv q(x)$ is the optimal candidate for a reaction coordinate.

In order to test the convergence of the algorithm, one has two different methods for estimating the probability
by other means:
\begin{itemize}
\item{} Perform a direct numerical simulation (DNS), that is, a direct Monte-Carlo approach. It is
done by simulating $N$ trajectories starting from ${\cal C}$
and monitor in one step those having reached ${\cal B}$ before ${\cal A}$.
Performing a large number of realisations yields an estimate of the crossing probability. This approach can be
achieved provided the temperature $1/\beta$ is not too small.
\item{} Solving (\ref{comm}) numerically. This is achieved by finite differences on a sufficiently large domain by
imposing some external Neumann boundary condition $\nabla q \cdot {\bf n} = 0$ on the boundary of the domain.
This is equivalent to forbid exits from the domain \cite{gar}. These boundaries are sufficiently far so that crossing
them is extremely unlikely even for reactive trajectories.
This approach is preferred over the Monte-Carlo
method when the temperature becomes too small. Note that, in 1-D, systems (\ref{eqpr}) always
derive from a potential and the committor has an explicit formula:
\begin{equation}
q(x)=\frac{{\displaystyle \int_{x_{\mathcal{A}}}^{x}\exp(\beta V){\rm d}x}}
{{\displaystyle \int_{x_{\mathcal{A}}}^{x_{\mathcal{B}}}\exp(\beta V){\rm d}x}}.
\label{comm1d}
\end{equation}
It can be further approximated using a saddle point approximation of the integral when $\beta$ is large by the formula:
\begin{equation}\label{comm_app}
 q(x) \simeq \frac{\sqrt{1-e^{-\omega (x_{\mathcal{A}}-x_s)^2}}+\Sigma\sqrt{1-e^{-\omega (x-x_s)^2}}}{\sqrt{1-e^{-\omega (x_{\mathcal{A}}-x_s)^2}}+\sqrt{1-e^{-\omega (x_{\mathcal{B}}-x_s)^2}}},~
\omega = -\frac12 \beta V''(x_s),
\end{equation}
where $x_s$ is the position of the saddle, and $\Sigma$ is the sign of $x-x_s$.
\end{itemize}
We now present the test models.
\begin{figure}
\centerline{\includegraphics[width=5cm]{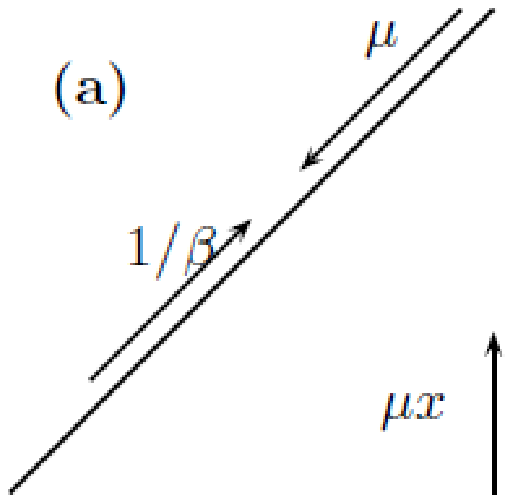}
\includegraphics[width=5cm,clip]{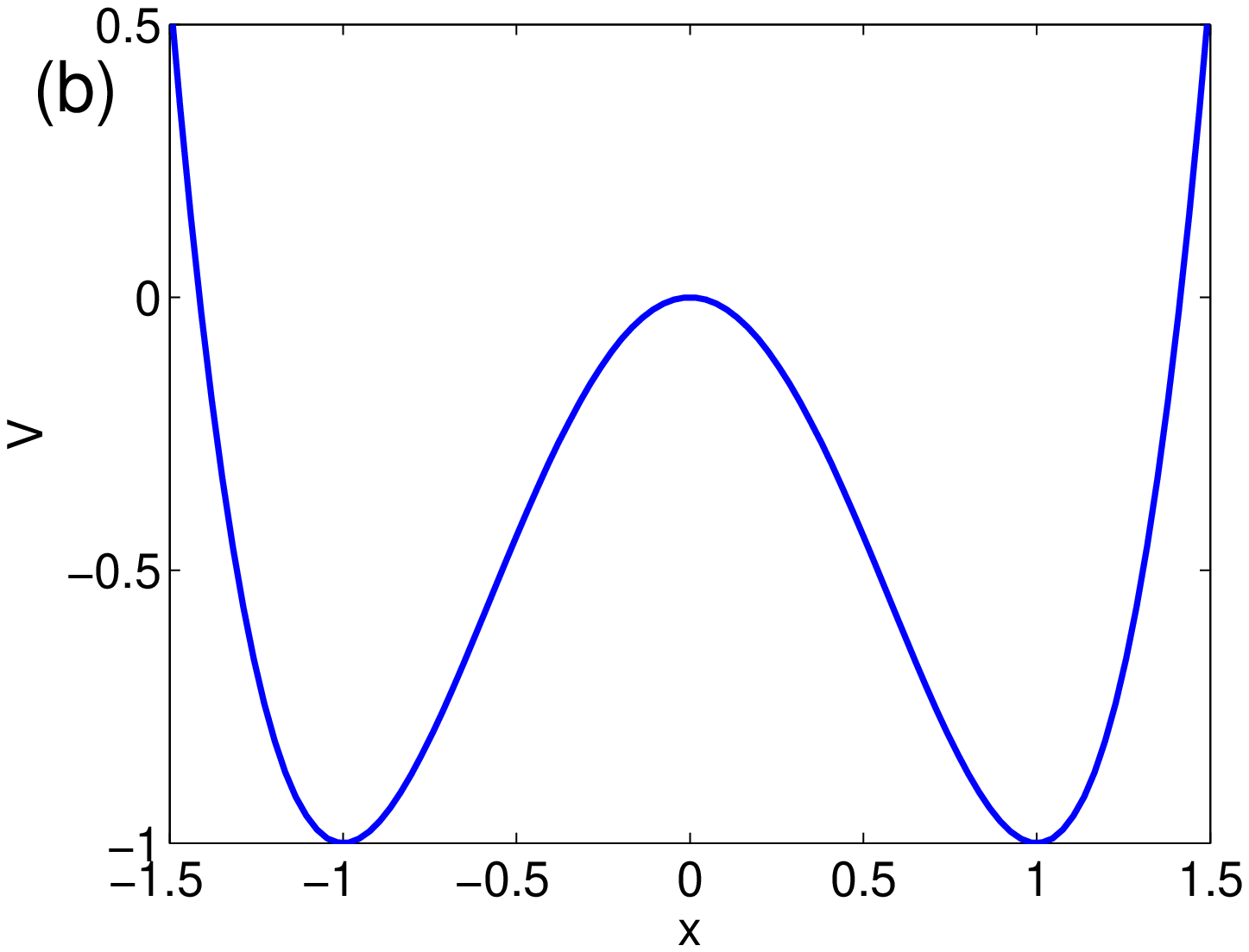}
\includegraphics[width=5cm]{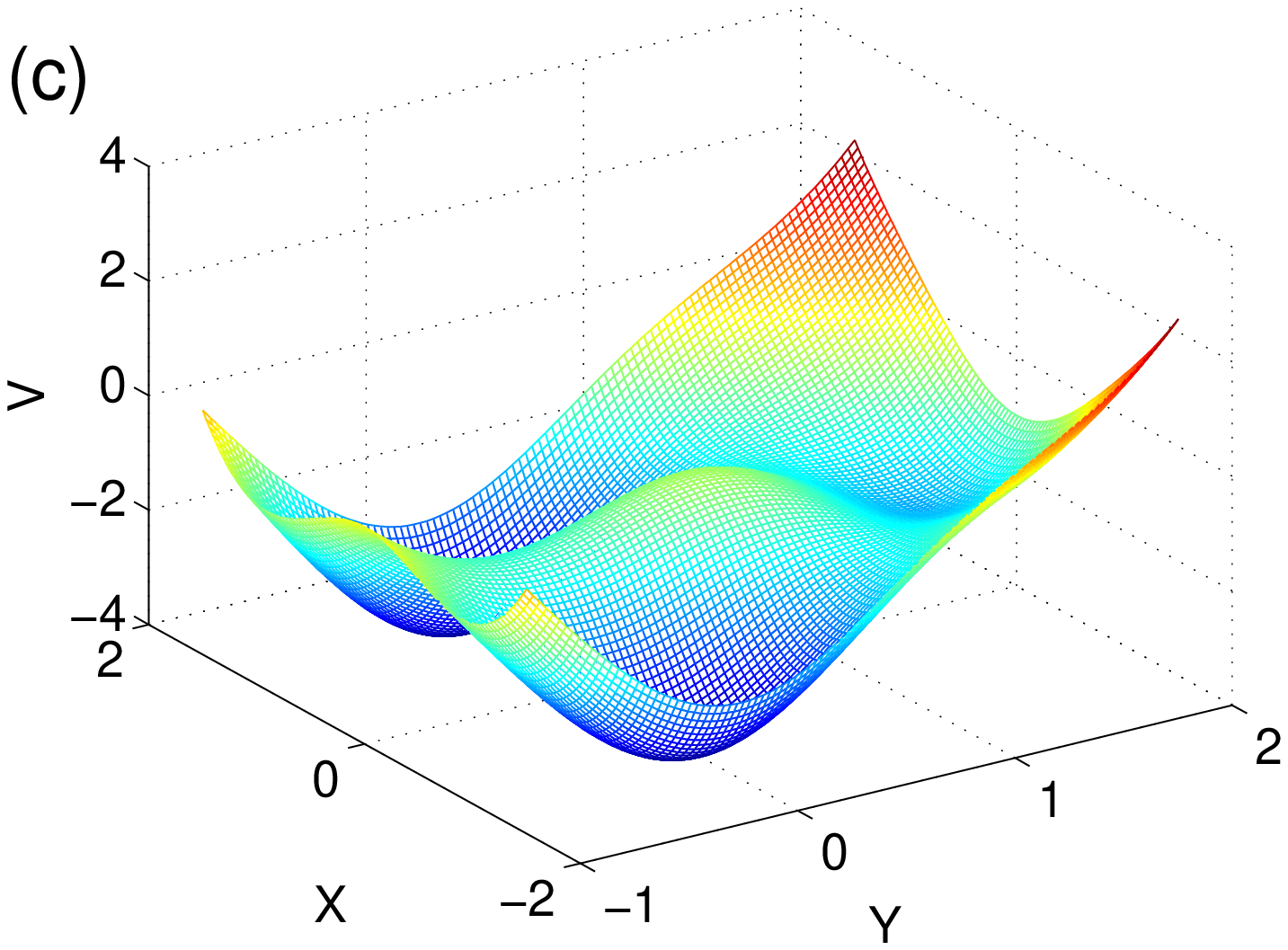}
\includegraphics[width=5cm]{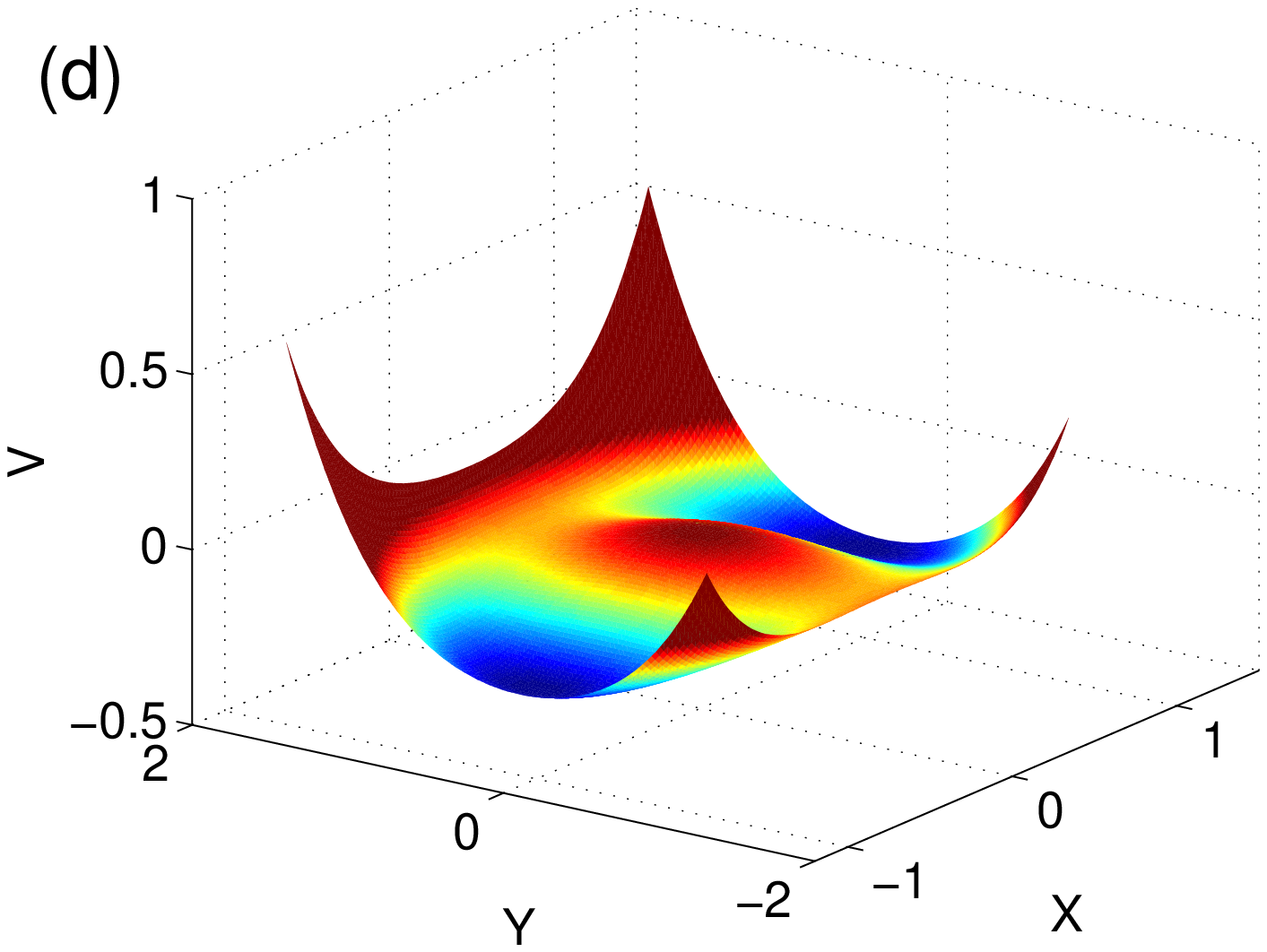}}
\caption{Sketch of the three models. (a): Brownian drift. (b): Potential of
the double well as a function of $x$. (c): colour levels
and surface plot of the potential of the triple well. (d):
colour levels and surface plot of the potential of the two saddles model.}
\label{fig0}
\end{figure}

\subsubsection{Brownian drift \label{brdr}}
The simplest non-trivial Langevin dynamics one can consider is the so-called Brownian drift \cite{cg07},
noted ``{\bf d}'' hereafter. In practice, it corresponds to a random walk in a moving frame, or by changing
the frame of reference, a Brownian motion with a constant drift force $-\mu$ (Fig.~\ref{fig0} (a)):
\begin{equation}
\frac{dx}{dt}=-\mu+\sqrt{\frac{2}{\beta}}\eta\,.
\label{drift}
\end{equation}
In this case, the sets $\mathcal{A}$ and $\mathcal{B}$, correspond
to positions $x_{\mathcal{A}}<x_0<x_{\mathcal{B}}$, and are not metastable.
The probability $\alpha$ is given by the committor function \cite{cg07} (Eq.~(\ref{comm1d})):
\begin{equation}
\alpha=q(x_0)=\frac{\sinh\left(\frac{\beta \mu}{2}(x_{\mathcal{A}}-x_0)\right)}{\sinh\left(\frac{\beta \mu}{2}(x_{\mathcal{A}}-x_{\mathcal{B}})\right)}\exp\left(\frac{\beta \mu}{2}(x_0-x_{\mathcal{B}}) \right)\,.
\label{com_drift}
\end{equation}
Provided $\mu>0$, reactive trajectories have a small probability $\alpha$ to occur.
Denoting $X_t^{x_0}$ the solution of (\ref{drift}) at time $t$ starting from initial condition $x_0$,
in dimension 1, $\alpha$ can be written as
\begin{equation}
\alpha = \Pr(X_\tau^{x_0} \geq x_{\cal B}),
\end{equation}
where $\tau = \min \{ t > 0,  X_t^{x_0} \notin ]x_{\cal A},x_{\cal B}[ \}$.

Simple reaction coordinates can be used: either the linear one, or a coordinate based
on the committor (Eq.~(\ref{com_drift})):
\begin{equation}
\Phi(x)=\frac{x-x_\mathcal{A}}{x_\mathcal{B}-x_\mathcal{A}}~{\rm or}~ \Phi(x)=q(x)\,. \label{coord_d}
\end{equation}
We choose $x_{\cal A}=0$, $x_0=1$, $x_{\cal B}=2$ and $\beta=1$ and vary $\mu$.
\subsubsection{1-D double-well potential}
The next level of complexity in term of rare events is the crossing of a saddle point in 1-D \cite{cglp}.
For that matter, we choose the dynamics of a particle in a symmetric double-well potential (Fig.~\ref{fig0} (b))
which is noted ``{\bf l1d}'' hereafter:
\begin{equation}
\frac{dx}{dt}=-\frac{dV}{dx}+\sqrt{\frac{2}{\beta}}\eta
\,,\, V(x)=x^4-2x^2\,.
\label{dw}
\end{equation}
There is only one control parameter: the
inverse of the temperature $\beta$. One has two metastable states $x=\pm1$, and one
saddle at $x=0$. The probability of crossing $\alpha$ can be computed semi-analytically using the
expression~(\ref{comm1d}), and is approximated by Eq.~(\ref{comm_app}).
We choose $x_\mathcal{A}=-1$, $x_\mathcal{B}=1$ and $x_\mathcal{C}=-0.9$.

As in the case of the drift, the choice of the reaction coordinate can be linear or the committor:
\begin{equation}
\Phi(x)=\frac{x+1}{2}~{\rm or}~ \Phi(x)=q(x) \,.\label{coord_1c}
\end{equation}
The approximation of the committor (Eq.~(\ref{comm_app})) can be used
for a simpler computation of $\Phi$. It can be shown that in 1-D, the choice of the reaction coordinate has no impact on the convergence of the crossing probability estimate. However, we will show that the convergence of the average duration trajectories does depend on the chosen reaction coordinate.

\subsubsection{2-D triple-well potential \label{tws}}
An other example which displays many interesting features is a 2-D
Langevin dynamics in a potential which has two symmetric global minima and one local
minimum \cite{cglp,VE} (referred to as ``{\bf l2d}'' in the following, see Fig.~\ref{fig0} (c)):
\begin{equation}
\frac{dx}{dt}=-\frac{\partial V}{\partial x}+\sqrt{\frac{2}{\beta}}\eta_x \,,\,\frac{dy}{dt}=-\frac{\partial V}{\partial y}+\sqrt{\frac{2}{\beta}}\eta_y \,,
\end{equation}
with:
\begin{equation}
 V(x,y)=0.2x^4+0.2\left(y-\frac{1}{3} \right)^2+ 3e^{\left(-x^2-\left(y-\frac{1}{3}\right)^2\right)}
 -3e^{ \left(-x^2-\left(y-\frac{5}{3}\right)^2\right)}  -5e^{ \left(-\left(x-1\right)^2-y^2\right)}-5e^{\left(-\left(x+1\right)^2-y^2\right)}\,.
\label{tw}
\end{equation}
Again, the behaviour of the system is controlled by the inverse of the temperature $\beta$.
We identify the two global minima at $(\pm1,0)$ with $\mathcal{A}$ and $\mathcal{B}$. These two sets will be defined quantitatively using the reaction coordinate, we set:
\begin{equation}
\mathcal{A}= \Phi^{-1}(-\infty,0.05]\,,\,\mathcal{B}=
\Phi^{-1}[0.95,+\infty)\,.\label{ab}
\end{equation}

The line $\mathcal{C}$ corresponds to $x=-0.9$ or $\sqrt{(x+1)^2+0.5y^2}=0.1$,
depending on our reaction coordinate.

One of the interesting feature in this model is the occurrence of a phase transition
yielding two distinct types of reactive trajectories as shown by the flux of reactive trajectories.

At high temperature $\beta < 5$, reactive paths
go through the saddle in the lower part of the domain between $\mathcal{A}$ and $\mathcal{B}$ at $(0,-0.3)$.
These paths have higher energy ($\Delta V=2.6$) but are shorter. The duration of trajectories
grows logarithmically with $\beta$ \cite{tps}.

At low temperature $\beta > 5$, the system is
controlled by energy rather than length. Most of the reactive trajectories will go through the
metastable minimum $\mathcal{D}$ at $(0,1.5)$ in the upper part.
They pass first through the intermediate saddle above ${\cal A}$ at
$(-0.6,1.1)$ (Fig.~\ref{fig0} (c)) and where
$\Delta V = V_{{\rm saddle}}-V_{\cal A} \approx 2.32$ to be trapped in the local minimum $\mathcal{D}$.
Then they go through the other saddle near ${\cal B}$ where
$\Delta V=V_{\rm saddle}-V_{\cal D} \approx 0.52$.
The trajectory has to escape the local minimum before reaching $\mathcal{B}$.
As a consequence, the trajectories include a first passage out of $\mathcal{D}$ and the duration of
trajectories increases this time exponentially with $\beta$.
Note that, The transition between these two regimes occurs smoothly at $\beta \simeq 5$.

\begin{figure}
\centerline{\includegraphics[width=6.5cm]{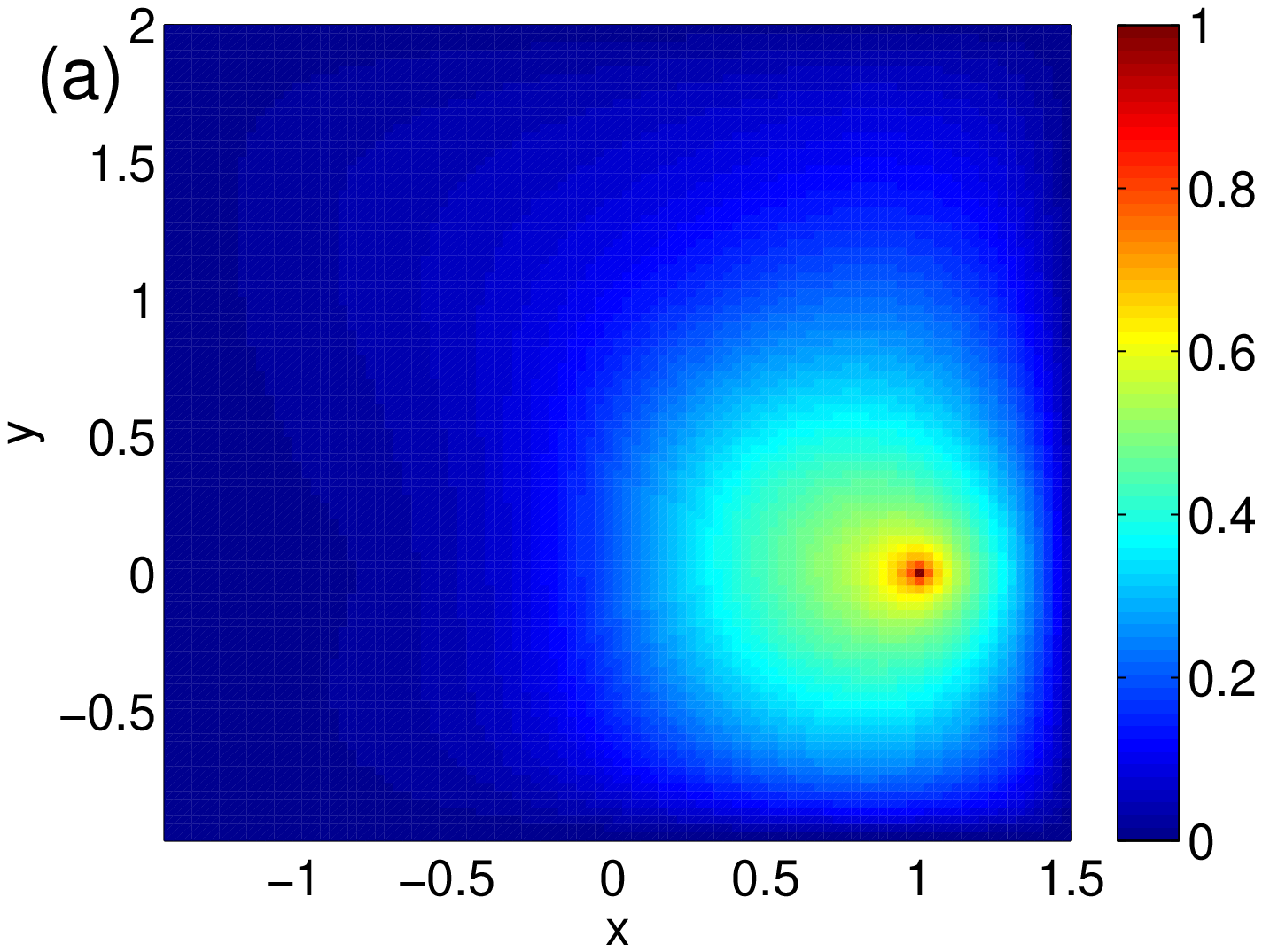}
\includegraphics[width=6.5cm]{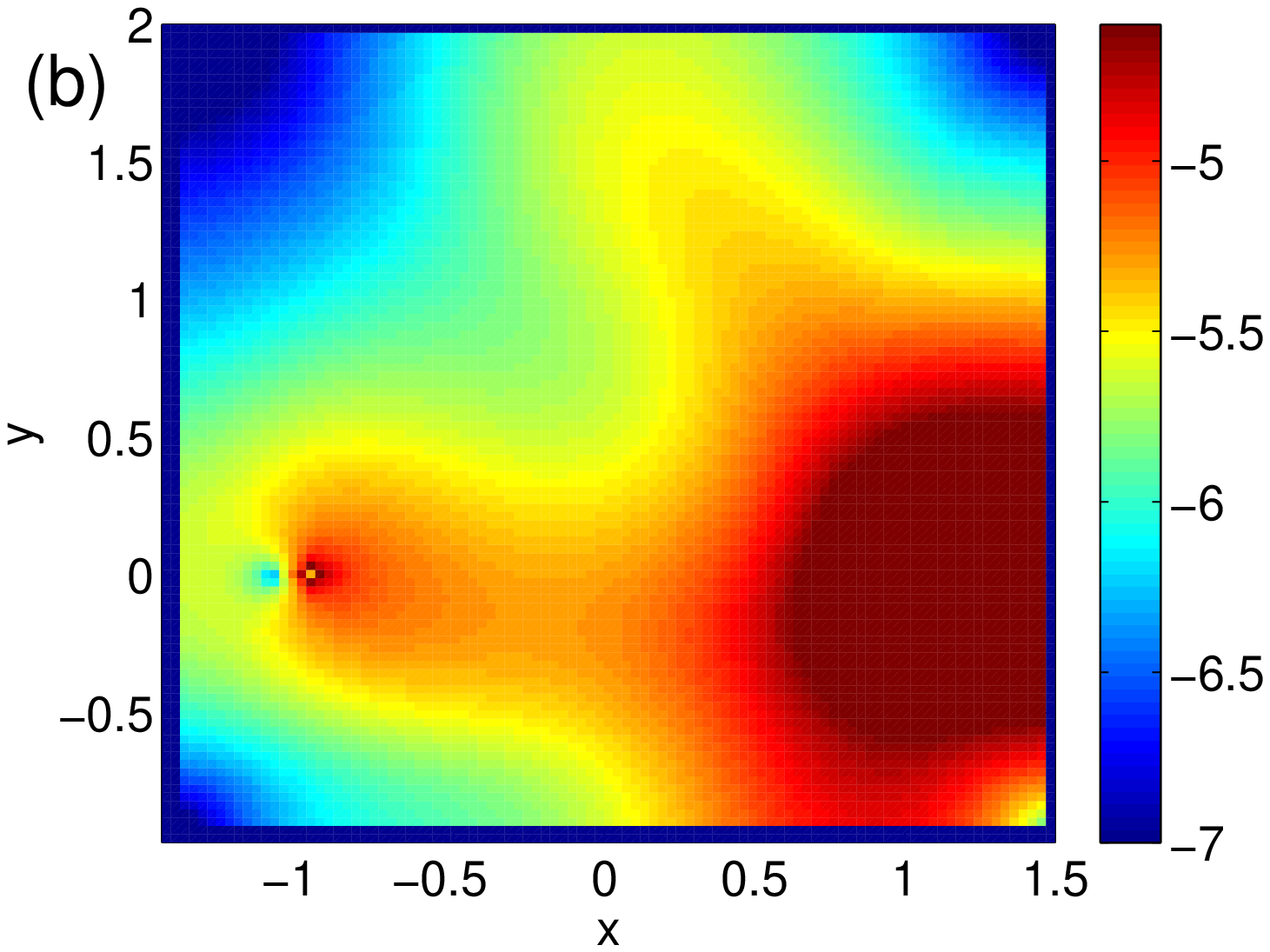}}
\centerline{\includegraphics[width=6.5cm]{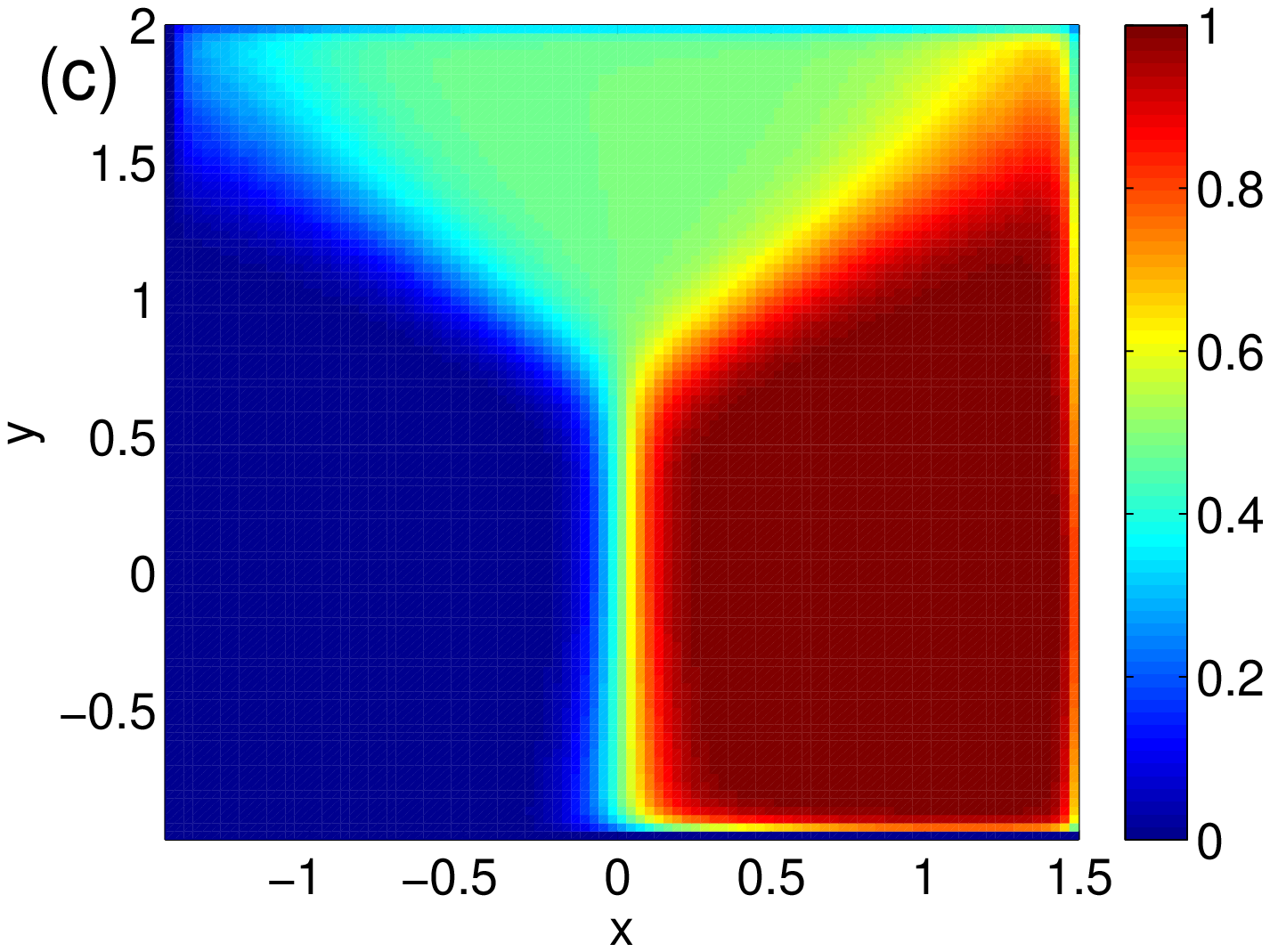}
\includegraphics[width=6.5cm]{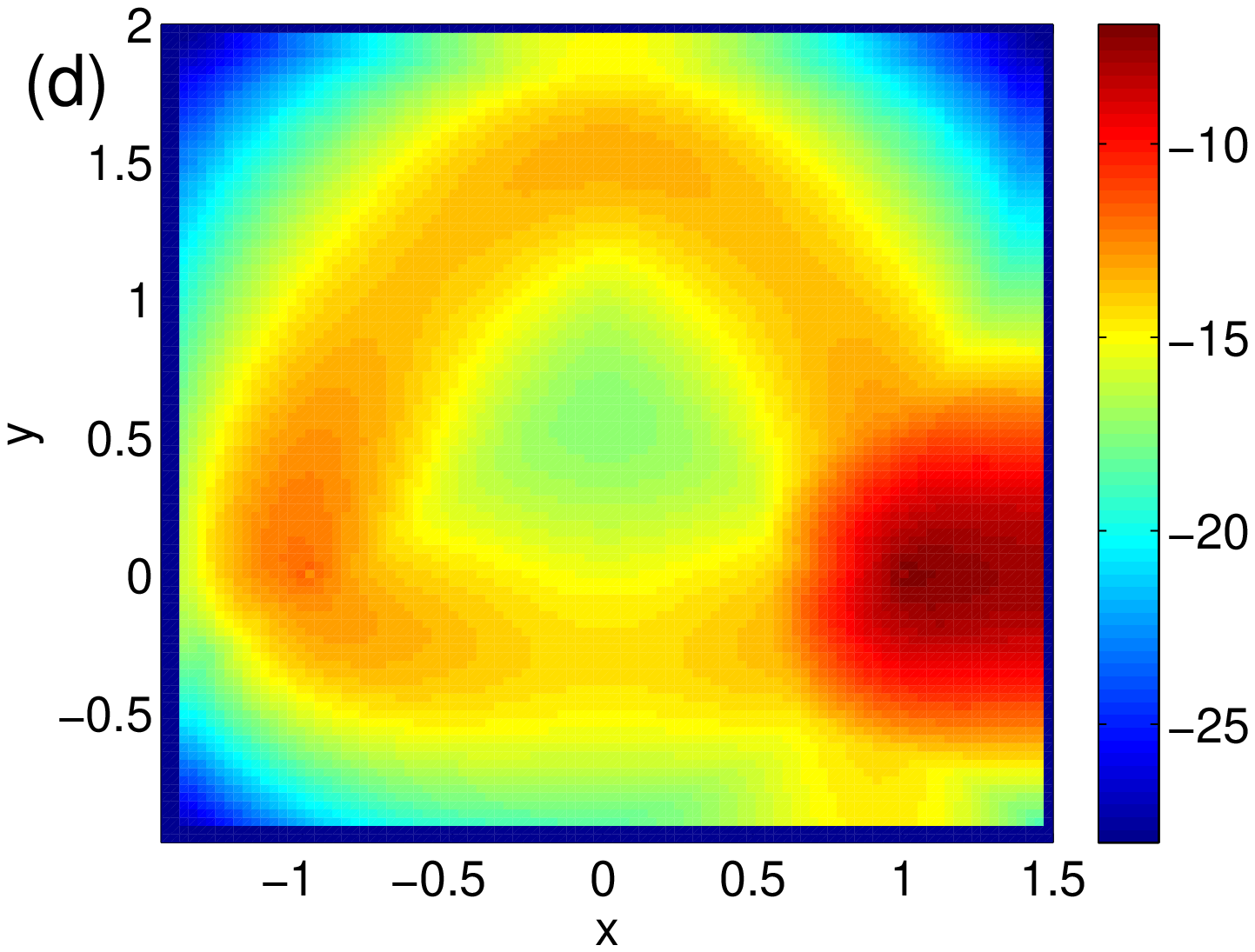}}
\caption{Committor function for the triple well model where ${\cal A}=\{ (-1,0)\}$
and $\mathcal{B}=\{ (1,0)\}$ for several value of the inverse temperature, and corresponding logarithm of the norm
of the flux of reactive trajectories. $\beta=1$: (a): committor, (b): norm of the flux. $\beta=10$:
(c): committor, (d): norm of the flux.}
\label{comm2d}
\end{figure}

Since the phase space is two-dimensional, the committor function can
also be computed directly by solving the backward Fokker-Planck equation numerically.
In our case, we solve equation (\ref{comm}) in the domain $[-1.5; 1.5]\times [-1;2]$.
We set $\mathcal{A}=(-1,0)$ and $\mathcal{B}=(+1,0)$ and the spatial mesh used is $\Delta x=\Delta y=0.03$.
The results are very similar to those found in \cite{VE} and are
shown in Fig.~\ref{comm2d}.
The sensitivity of the committor to the definition
of $\partial\mathcal{A}$ and $\partial\mathcal{B}$ decreases rapidly as $\beta$ is increased.
The computation at $\beta=1$ (Fig.~\ref{comm2d} (a)) shows the near ``free diffusion''
character of the trajectories at high temperature. Since temperature is relatively high,
the system does not really feel the potential $V$
and the committor is nearly that of a random walk between two points of the plane.
The two possible paths between $\mathcal{A}$ and $\mathcal{B}$
appear more and more clearly as $\beta$ is increased. The lower
path is characterised by the steep increase of the probability of crossing as one goes through the lower saddle,
while the upper channel is characterised by a flat plateau of $q\simeq 1/2$ in the
neighbourhood of the metastable minimum (Fig.~\ref{comm2d} (c)).
Computation of the flux of reactive trajectories (Fig.~\ref{comm2d} (b,d)) is in agreement with
previous computations and clearly shows the change of type of trajectories as $\beta$ is increased \cite{VE}.

Since the model is two-dimensional, there are many possible reaction coordinates.
One can choose a simple linear coordinate, or a coordinate weighted by $y$, based on the Euclidian norm \cite{cglp}:
\begin{equation}
\Phi_l(x)=\frac{x+1}{2}\,,\, \Phi_n(x)=\frac{1}{2}\sqrt{(x+1)^2+\frac{1}{2}y^2}\,.\label{coord_2c}
\end{equation}
This second reaction coordinate (hereafter denoted as ``norm'') is a simple way of mimicking
the behaviour of the committor function,
since it gives some weight to trajectories on the upper half plane $y>0$.
Eventually, one can use a reaction coordinate $\Phi_c$ based on the committor itself.
In order to do that, one first solves equation (\ref{comm}) (see Fig.~\ref{comm2d}) and
records values of the committor function. During the
AMS algorithm simulation, an estimate of the committor at any position on the plane
is obtained by simple interpolation from the numerical solution of equation~(\ref{comm}) (Fig.~\ref{comm2d}). When using the committor as reaction coordinate, we use for consistency
the same definitions of $\mathcal{A}$, $\mathcal{B}$ and $\mathcal{C}$ given by the norm reaction coordinate.

\subsubsection{Two saddles}
  We consider a simpler type of 2-D Langevin dynamics using a potential that has two minima ($(x=\pm 1,y=
0)$), one local maximum $(x=0,y=0)$ and two saddles between the minima $(x=0,y=\pm1)$ (Fig.~\ref{fig0} (d)):
\begin{equation}
V=\frac{x^4}{4}-\frac{x^2}{2}-0.3\left( \frac{y^4}{4}-\frac{y^2}{2} +x^2y^2\right)\,,\,
\frac{dx}{dt}=-\frac{\partial V}{\partial x}+\sqrt{\frac{2}{\beta}}\eta_x \,,\,
\frac{dy}{dt}=-\frac{\partial V}{\partial y}+\sqrt{\frac{2}{\beta}}\eta_y\,.
\end{equation}
The crossing in this model is simple: it goes from one minimum to the other, \emph{via} one of the
two saddles, with a probability $\frac{1}{2}$ to go through each channel due to the symmetry $y \to - y$.
Its main interest will be
in the discussion of the quality of the reaction coordinate in the triple-well potential model. Indeed, it
will help infer more generic behavior for the variance.
The reaction coordinates used are:
\begin{equation}
\Phi_l(x)=\frac{x+1}{2}\,,\, \Phi_n(x)=\frac{1}{2}\sqrt{(x+1)^2+\frac{1}{2}y^2}\,.\label{coord_2s}
\end{equation}
We take the same definition of $\mathcal{A}$ and $\mathcal{B}$ as in the triple well case (Eq.~(\ref{ab})) and the line
$\mathcal{C}$ corresponds to $x=-0.9$ or $\sqrt{(x+1)^2+0.5y^2}=0.1$, depending on the choice of reaction coordinate
(respectively $\phi_l$ or $\phi_n$).

\subsection{Algorithm Complexity}

The complexity $C$ of the algorithm reads:
\begin{equation} C=\underbrace{N(D_{init}+S_1)}_{stage 1})+
\underbrace{k(nD_{branch}+S_2)}_{stage 2}\,, \label{comp}\end{equation}
with $NS_1$ the complexity of the initial sorting and
$S_2$ that of replacement/sorting at each step. The cost for simulating a single trajectory
either during the initial step or during the branching step is represented
by $D_{init}$ and $D_{branch}$ respectively.
The use of an optimal sorting/replacement algorithm yields $S_1,S_2=O(\ln N)$.
If $n/N$ is small and by approximating $\ln(\alpha)\propto -\beta$, we find
using Eq.~(\ref{alpha_algo}):
\begin{equation}K=\ln(\alpha N/r)/\ln(1-n/N)\simeq (N/n)|\ln(\alpha N/r)|.\end{equation}
This yields:
\begin{equation} \frac{C}{N}=(D_{init}+a \ln(N))+
\frac{1}{n}|\ln(\frac{\alpha N}{r})|(nD_{branch}+b \ln(N))\,. \label{compf}\end{equation}
In practice $N/r \simeq 1$, that is $K \propto N\beta/n$.

The costs $D_{init}$, $D_{branch}$ strongly depend on the model, \emph{i.e.}
the number of degrees of freedom, duration of trajectories, and numerical precision one wants to achieve. It
can be large or small relatively to $S_1$ and $S_2$.
If it is large, the complexity is relatively insensitive to the number of suppressed trajectories $n$.
In the case of the drift and double well models,
the complexity $C$ increases only with the number of clones:
since the numerical simulations are rather short, the algorithm is very sensitive to the
complexity of the sorting and replacement procedures. In the case of the triple well model,
the complexity $C$ increases nearly linearly with the number of clones $N$, \emph{via} the number of iterations,
and the maximum duration of reactive trajectories.
It is relatively insensitive to the cost induced by sorting the maxima of trajectories.

We can consider the strategy of suppressing a proportion $n/N$ of trajectories. This approach is particularly
interesting when applying the algorithm to stochastic partial differential equations (SPDEs)
with many degree of freedoms. In
this case, one can parallelise not only the first step of the algorithm, by distributing the $N$ clone dynamics,
but the branching of the trajectories as well, by distributing the $n$ evolutions. It allows a possible reduction of the complexity of
the algorithm (Eq.~(\ref{comp})), by replacing $D_{init,branch}$ by $D_{init,branch}/\sharp$,
with $\sharp\le n$ the number of cores available.

\subsection{Time discretisation \label{dis}}
We briefly discuss the numerical approximation of the stochastic equation needed to generate trajectories.
One can use a simple Euler approximation
\cite{gar}, or high-order schemes \cite{KP}.
High-order schemes involve a more complicated treatment of the noise term \cite{KP,gar}.
In practice, if we define two normal random variables $U_{1,2}$ and two discretisations of the Wiener processes $\eta$: $\Delta W$, $\Delta Z$. They are defined as:
\begin{align} \notag
U_{1,2}\sim N(0,1)\,,\, \Delta W=\sqrt{\frac{2dt}{\beta}}U_1\,,\, \Delta Z=\sqrt{\frac{2}{\beta}}\frac{1}{2}(dt)^{\frac{3}{2}}\left(U_1+\frac{1}{\sqrt{3}}U_2 \right)\,,
\end{align}
such that
\begin{equation}
E(\Delta Z)=0\,,\,E(\Delta Z^2)=\frac{2}{3\beta}dt^3\,,\,E(\Delta W\Delta Z)=\frac{1}{\beta}dt^2\,.
\end{equation}
The Euler approximation of equation (\ref{eqpr}) reads:
\begin{equation}
X_{n+1}=X_n+F(X_n)dt+ \Delta W \label{05}
\end{equation}
and is valid for $X$ being either a scalar or a vector. The strong order of convergence of this numerical scheme is $dt^{\frac{1}{2}}$.
One can use a scheme whose order is $dt^{\frac{3}{2}}$ if $X\in \mathbb{R}$ :
\begin{equation}
X_{n+1}=X_n+F(X_n)dt+\Delta W+F'\Delta Z+\frac{1}{2}dt^2(FF'+\frac{1}{\beta}F'')\,. \label{15}
\end{equation}
In practice, the Euler scheme~(\ref{05}) is used for all four models, while the high-order scheme
(\ref{15}) is used for the double-well potential. Note that the rate of convergence of theses schemes does
not presume of the rate of convergence of an estimate of $\alpha$ as a function of $dt$.
Indeed, they are defined as the rate of convergence of the ensemble average of the error.
They do not give information on the hitting times $t_{\mathcal{A},\mathcal{B}}$ of fluctuations, on which the algorithm is based. In practice, one can show that these hitting times converge with $dt^{\frac{1}{2}}$ \cite{ght} with any kind of Euler Scheme. One should not expect a better order for $\alpha$.

\section{Numerical convergence \label{S2}}
In this section, we study the convergence of the estimator $\hat \alpha$ given in (\ref{alpha_algo})
as a function of the timestep $dt$ and
the number of particles $N$. We will denote the estimator $\hat{\alpha}=\hat{\alpha}_{dt,N}$.
We also consider the convergence of the average duration of reactive trajectories $\tau$.
We write the estimator provided by the algorithm as
\begin{equation}
\hat{\alpha}_{dt,N}=\alpha_{{\rm DNS}} + b + v, ~{\rm with}~\alpha_{\rm DNS} \equiv \alpha(1-e_{dt}),
\end{equation}
where $b$ stands for the bias with respect to the DNS value of $\alpha$, $v$ a centered
random variable accounting for the observed variance and $e_{dt}$ is the discretisation error.
We focus here on the term $e_{dt}$. For 1-D systems, The behaviour of $v$ (with zero mean) is well understood
and the bias can be corrected through reformulations of the algorithm \cite{com}.
The value of $\alpha$ is computed analytically, semi analytically (Eq.~(\ref{comm1d})) and/or by DNS
by letting $dt\rightarrow 0$ and $N\rightarrow \infty$.
We use a large range of parameters for the drift and the double-well potential: $\mu$ goes up
to $40$ ($\alpha\simeq4\cdot 10^{-18}$) and $\beta$ up to $20$ ($\alpha\simeq9\cdot 10^{-10}$) in the
computation of $\hat{\alpha}_{dt}$. In the case
of the triple-well potential, we consider the case of small $\beta$ and leave small probabilities to section~\ref{reac}.

\subsection{Time step}
The convergence with the time step $dt$ is quantified by the rate $\gamma$ at which the numerical error $e_{dt}$ goes to zero:
\[e_{dt} \equiv\left|1-\frac{\langle\hat{\alpha}_{dt}\rangle-b}{\alpha}\right| \approx dt^{\gamma}\,, \]
where $b$ is the bias of the estimation. For that matter, several realisations algorithm are performed with a
large value of $N=20000$ in a situation where both bias and variance are relatively small.
The crossing probability $\hat{\alpha}_{dt}$ is averaged over these realisations to reduce the size of the interval of confidence.

In order to verify that the bias can be neglected, we compute $\alpha_{{\rm DNS},dt}$, provided it is not too small. We
find that the DNS
has the same rate of convergence as the algorithm, and that the relative difference
$b=1-\frac{\hat{\alpha}_{dt}}{\hat{\alpha}_{{\rm DNS},dt}}$, of order $10^{-3}$, does not
depend on $dt$ in the regime considered.
\begin{figure}
\centerline{\includegraphics[height=5.2cm]{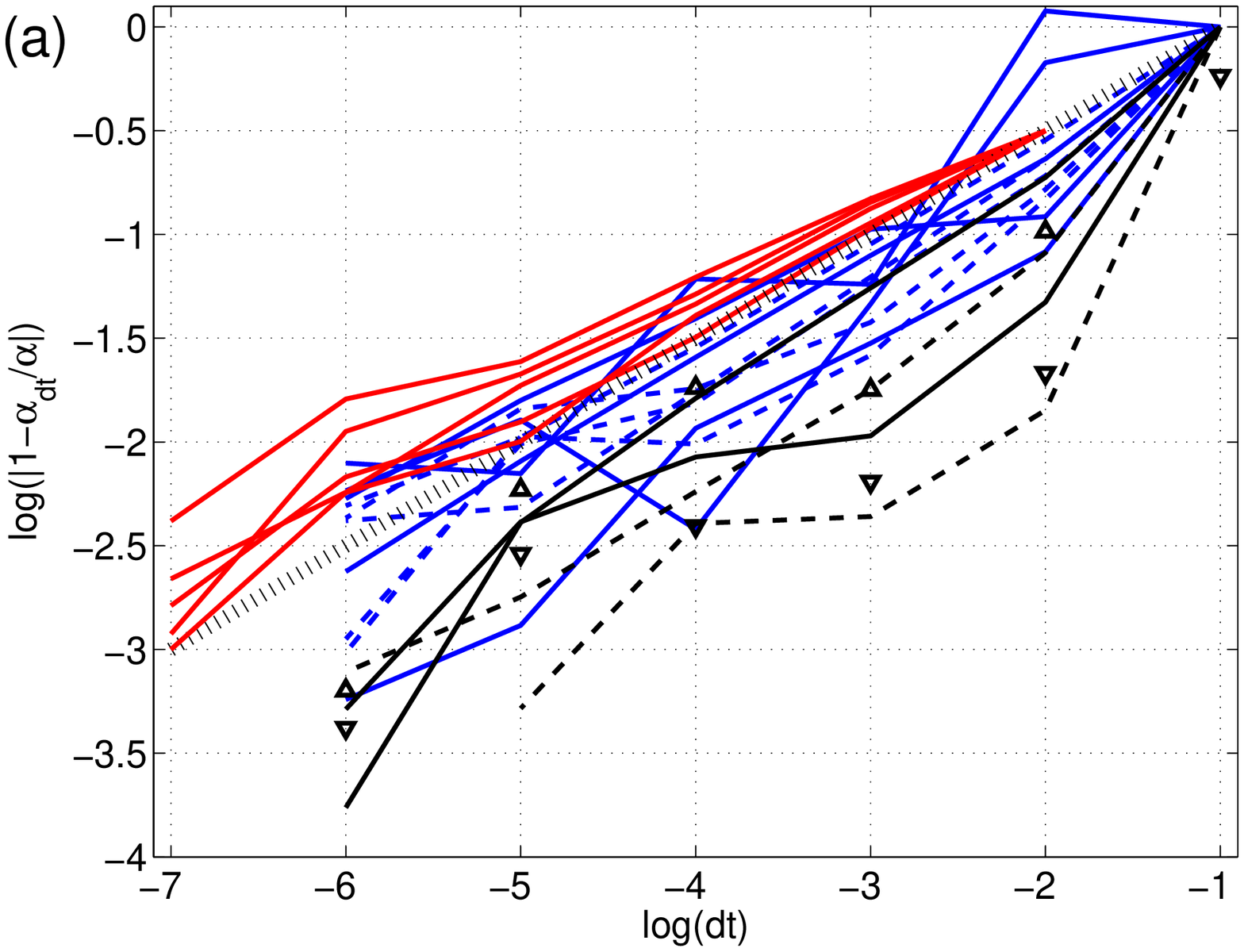}
\includegraphics[height=5.2cm]{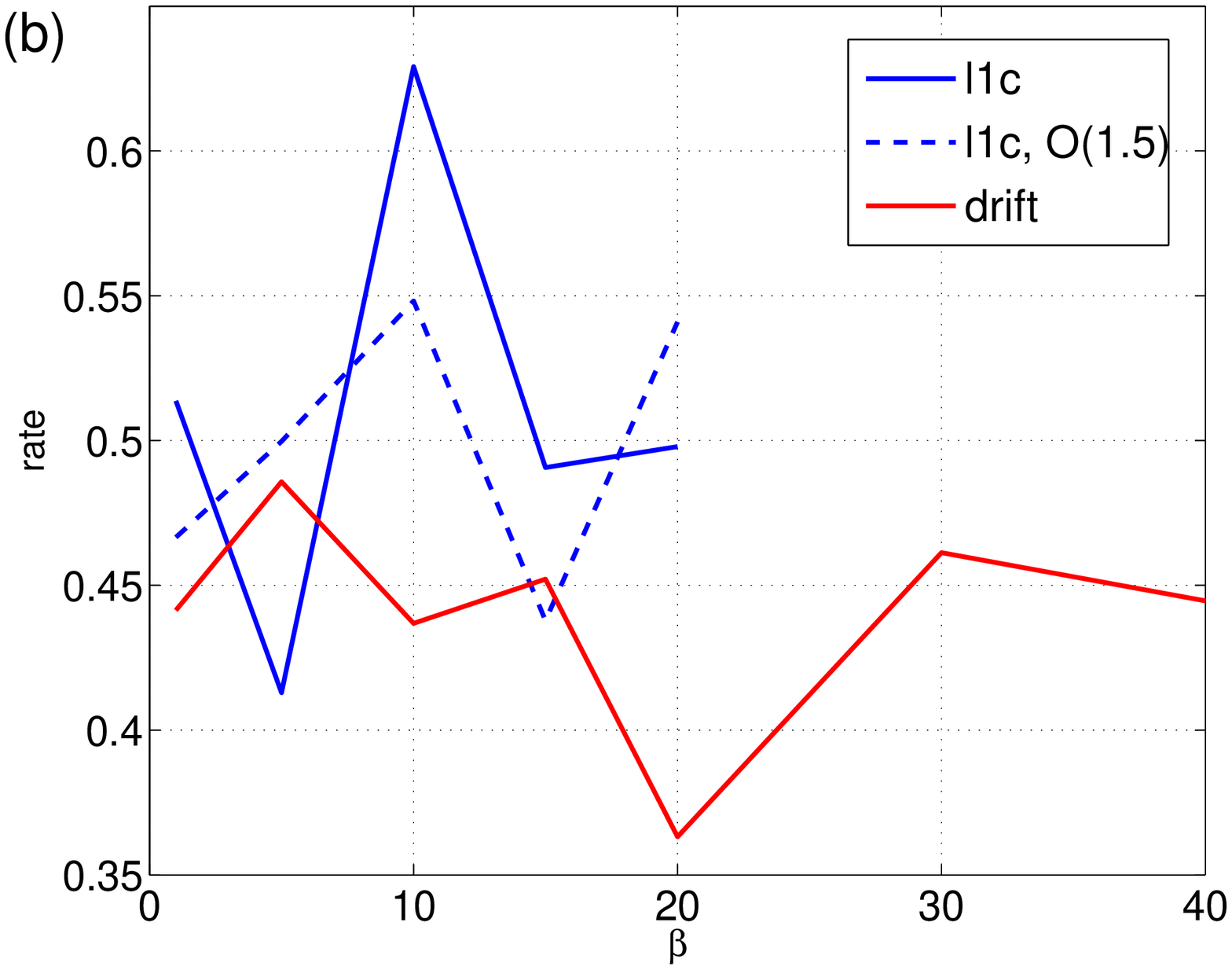}
\includegraphics[height=5.2cm]{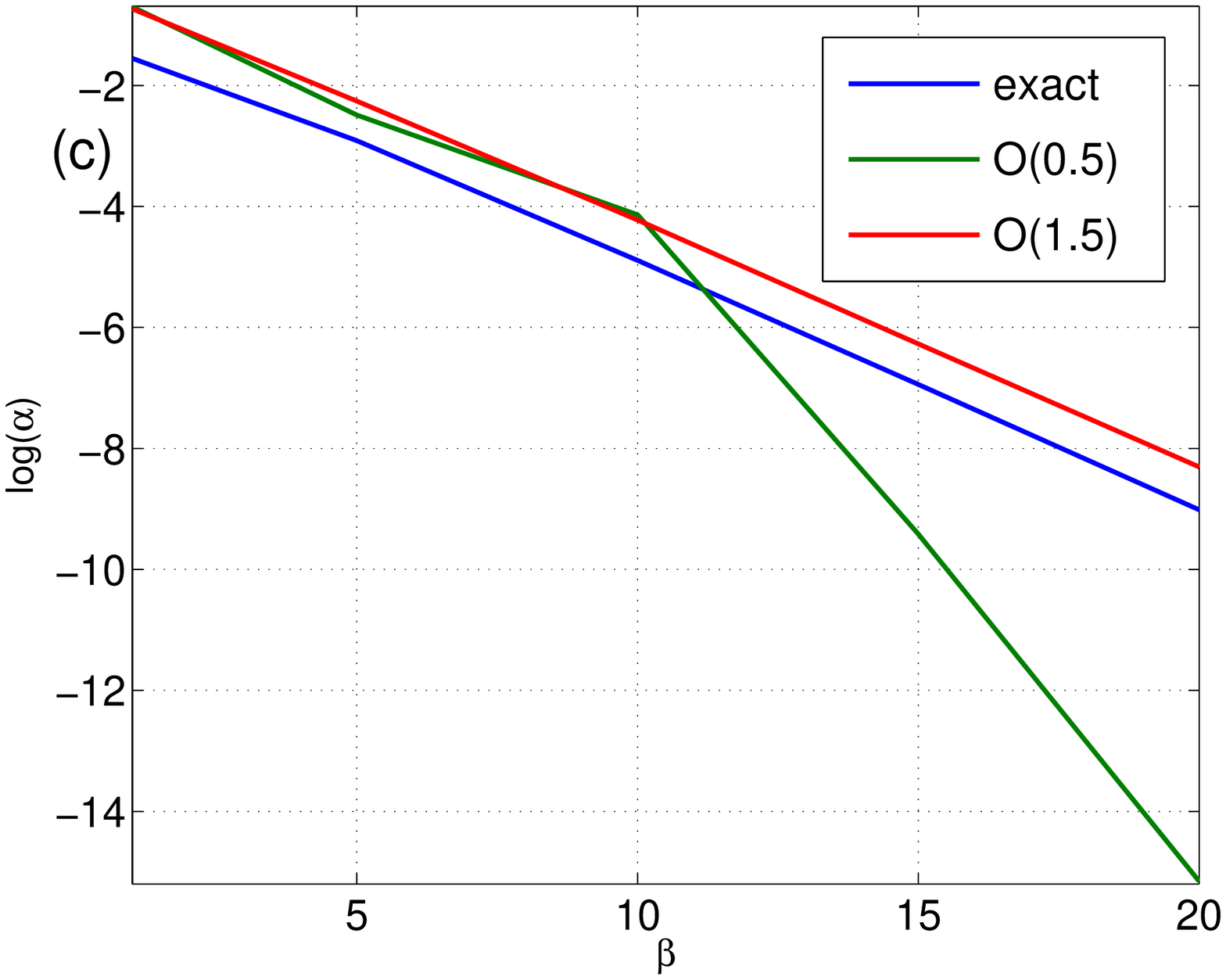}}
\caption{(a): convergence of the computed transition probability toward the analytical value, as a
function of $dt$ in $\log$ scale, for the double well, blue (two numerical schemes, Euler, full line, order $1.5$, dashed line), drift (red) and triple well
(black, full line: linear reaction coordinate, dashed line: norm reaction coordinate)
models for increasing values of $\beta$ or $\mu$. The error bars are obtained
from the confidence interval of $\alpha_{dt,N}$. (b): Rate of convergence of the probability with $dt$ as function of $\beta$ for
the double well model (two different numerical schemes) and of $\mu$ the drift.
(c): comparison between Euler and order 1.5 scheme as a function of $\beta$ for large $dt = 0.1$.
As $\beta$ increases Euler scheme strongly underestimates the probability.}
\label{rate}
\end{figure}

We display $\log(e_{dt})$ for all models and parameters as a function of $\log(dt)$ in Fig.~\ref{rate} (a).
The value of $\alpha$ is obtained either from analytics or when possible with DNS with a very small timestep.
The error bars are obtained from the variance of $\hat{\alpha}_{dt}$. They are placed when relevant. Indeed, for many values of $\beta$ or $dt$, the variance is two small for the bars to be visible; meanwhile, for $dt=0.1$, the variance is often extremely large and the bars are out of the figure. This last case is detailed in the next paragraph.
One finds that $\gamma\simeq 0.5$ in all the cases and numerical schemes.
The rate of convergence $\gamma$ is obtained by a fit of $\log(e_{dt})$ and is displayed in Fig.~\ref{rate} (b)
as a function of $\beta$ and $\mu$ for the double-well potential and the drift respectively.
It appears that $e_{dt}\propto \sqrt{dt}$, independently of the value
of $\alpha$ or the strong order of convergence of the numerical scheme. This shows that $\hat{\alpha}$
follows the same order of convergence as the hitting times \cite{ght}. We moreover find that $\alpha_{\rm DNS}$
follows the same power law, which confirms the fact that the error is not rooted in the algorithm itself.
Note that the discretisation error can lead
to a much larger error than the bias and the variance in the estimation of $\alpha$: one needs
$dt\simeq 10^{-6}$ to have a relative error of order $10^{-3}$. A small $dt\le 10^{-4}$ and
a large $N$ are both necessary for a precise estimate of $\alpha$.

In practice the effect of the order of convergence in the numerical scheme is seen when
using relatively large time steps.
For large $dt\simeq 0.1$ and large $\beta$ or $\mu$, the number of steps of the
algorithm is extremely high, which yields probabilities much smaller than
$\alpha$, as well as a very large variance (see the error bars at $dt=0.1$).
The consequence is a poor result and a very long
run time of the algorithm. This problem is solved by the use of the order $1.5$ scheme in the case of
the double-well potential (see Fig.~\ref{rate} (c)).

On a side note, we have also implemented a Brownian bridge for 1-D models and Euler scheme
in order to improve the branching procedure. In practice, when branching trajectory $i$ on trajectory $i^\star$, one has $X^{i^\star}_{t^\star -dt}<l_k<X^{i^\star}_{t^{\star}}$. The Brownian bridge is run
between $X^{i^\star}_{t^\star -dt}$ and $X^{i^\star}_{t^{\star}}$ with a time step $dt/100$.
It produces ${t^\star}'$  such that $X^{i^\star}_{{t^\star}' -dt/100}<l_k<X^{i^\star}_{{t^{\star}}'}$.
The new trajectory is then simulated starting from the initial condition $x_{i^\star}'=X^{i^\star}_{{t^{\star}}'}$.
It leads to no significant improvement (not shown) and therefore is not used in the following.

\subsection{Bias and variance on the estimation of  and $\tau$ \label{se_norm}}
\begin{figure}
\centerline{\includegraphics[width=6.5cm]{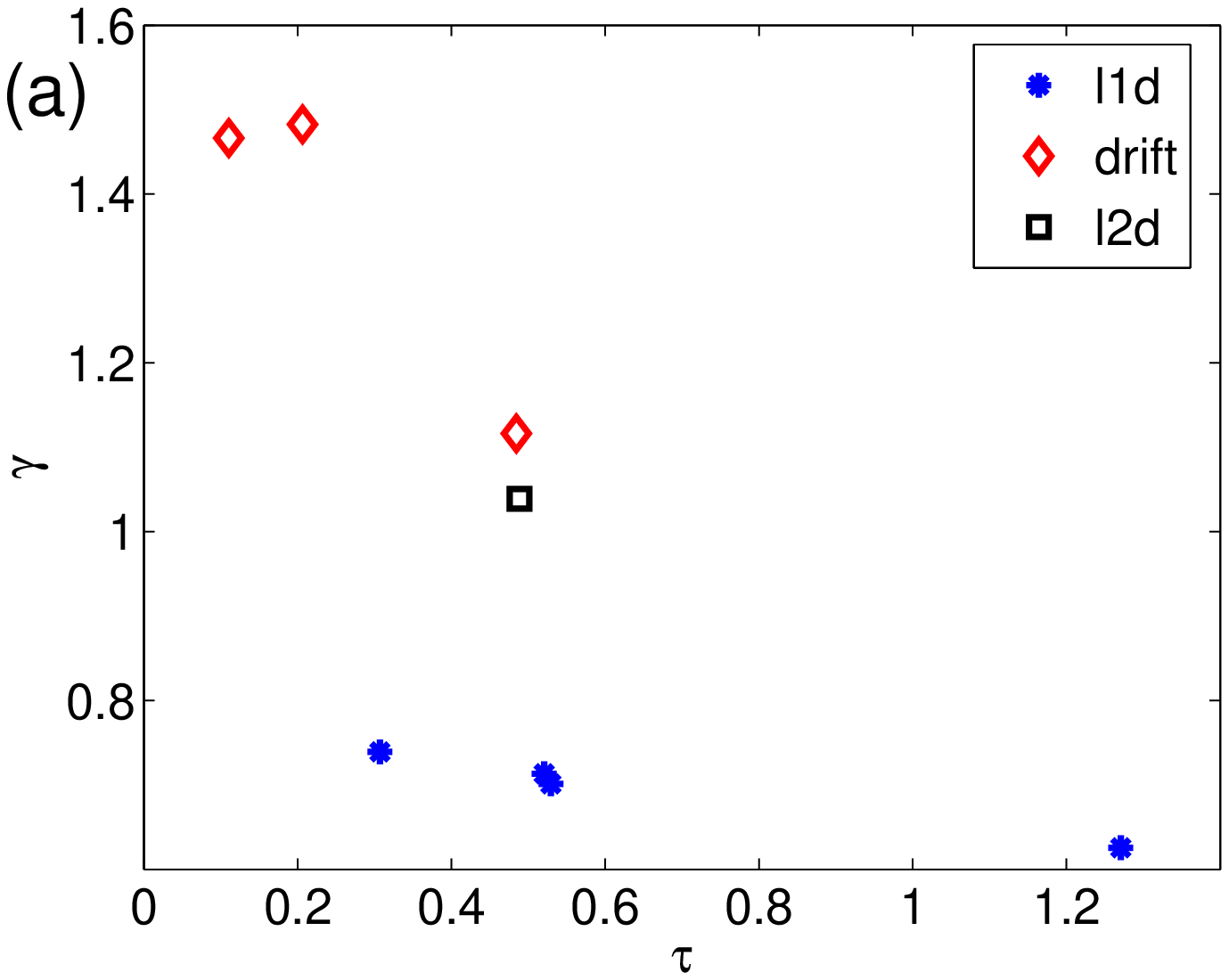}
\includegraphics[width=6.5cm]{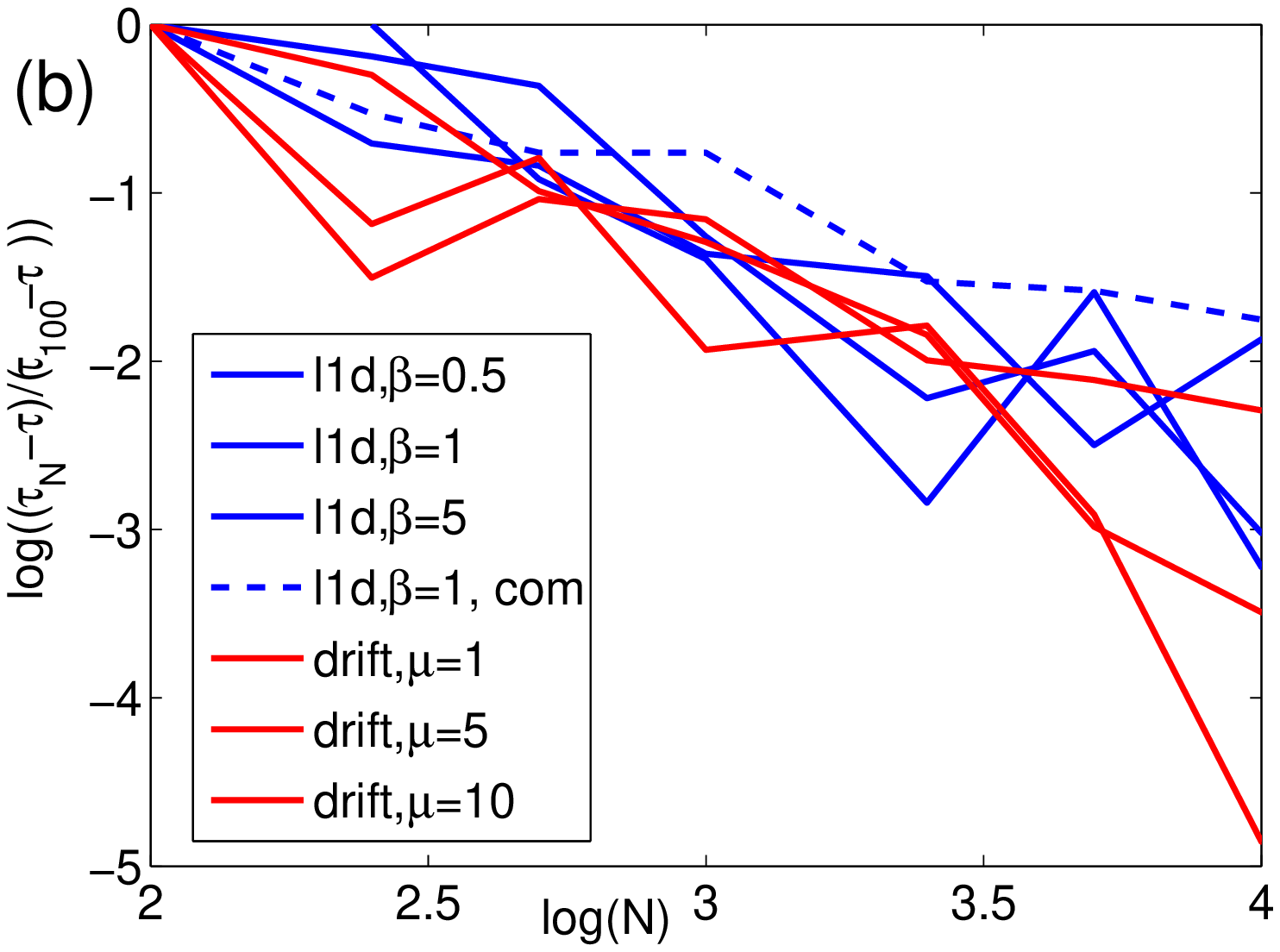}}
\caption{Convergence of the average duration of a reactive trajectory with the number of clones. (a) : $\gamma$ as a function of $\tau$. (b) : logarithm of the bias as a function of $\log(N)$ .}
\label{stat_temps}
\end{figure}

We now investigate estimates of the duration $\tau$ of reactive trajectories. This is
motivated by the fact that systematic theoretical predictions on this quantity are lacking.
Theoretical results mostly concern 1-D systems with continuous trajectories and in the limit of
large $\beta$ \cite{tps}. However, for rather general situations,
the distribution of $\tau$ appears to be of Gumbel type \cite{cg07,cglp,dur} (Fig.~\ref{triple_} (b)).
We use the same settings as for the statistics of the crossing probability. We compute the distribution of $\langle\tau\rangle_N$ for $n_c$ values of $100\le N\le 20000$ over at least $2000$ realisations of the algorithm. The variance of the durations appears to decrease with $N$. Since we expect a central limit behaviour for this quantity, it is compared to $\langle\tau_{DNS}\rangle/\sqrt{N}$, by computing:
\begin{equation}
\gamma=\frac{1}{n_c}\sum_{N}\frac{\hat{\sigma}_{\tau,N}\sqrt{N}}{\tau_{DNS}}
\end{equation}
as a function of $\tau$. One can see that there is a consistent behaviour for
each model, \emph{i.e.} in good agreement with the proposed dependence, $\gamma$ is always of order one,
with small differences between each model (Fig.~\ref{stat_temps} (a)). The average durations
are all negatively biased with respect to the average duration computed \emph{via} DNS. The logarithm
of the bias is displayed as a function of the logarithm of $N$ in Fig.~\ref{stat_temps} (b).
It decreases like $1/N^2$ for all cases except that of the double-well potential using the
committor reaction coordinate, which behaves as $\tau-\langle\hat{\tau}\rangle\propto N^{-3/2}$.

\section{Slow convergence and generic behaviour of the variance \label{reac}}
In this section, we focus on AMS behavior in dimension larger than one. We therefore study
the triple-well and two saddles models.
We distinguish two different situations. In the first case, we focus on the effect of
a phase transition on the algorithm convergence. In the second case, we rather consider
how the AMS estimate degrades in the limit of small noise.
A three-level simplification of the triple-well potential is proposed at the end. It helps
 understanding the effect of the modified
selection of trajectories on the estimation of the average duration and to quantify the rate of convergence.

\subsection{Ensemble of trajectories during phase transition}
\begin{figure}
\centerline{\includegraphics[height=5cm]{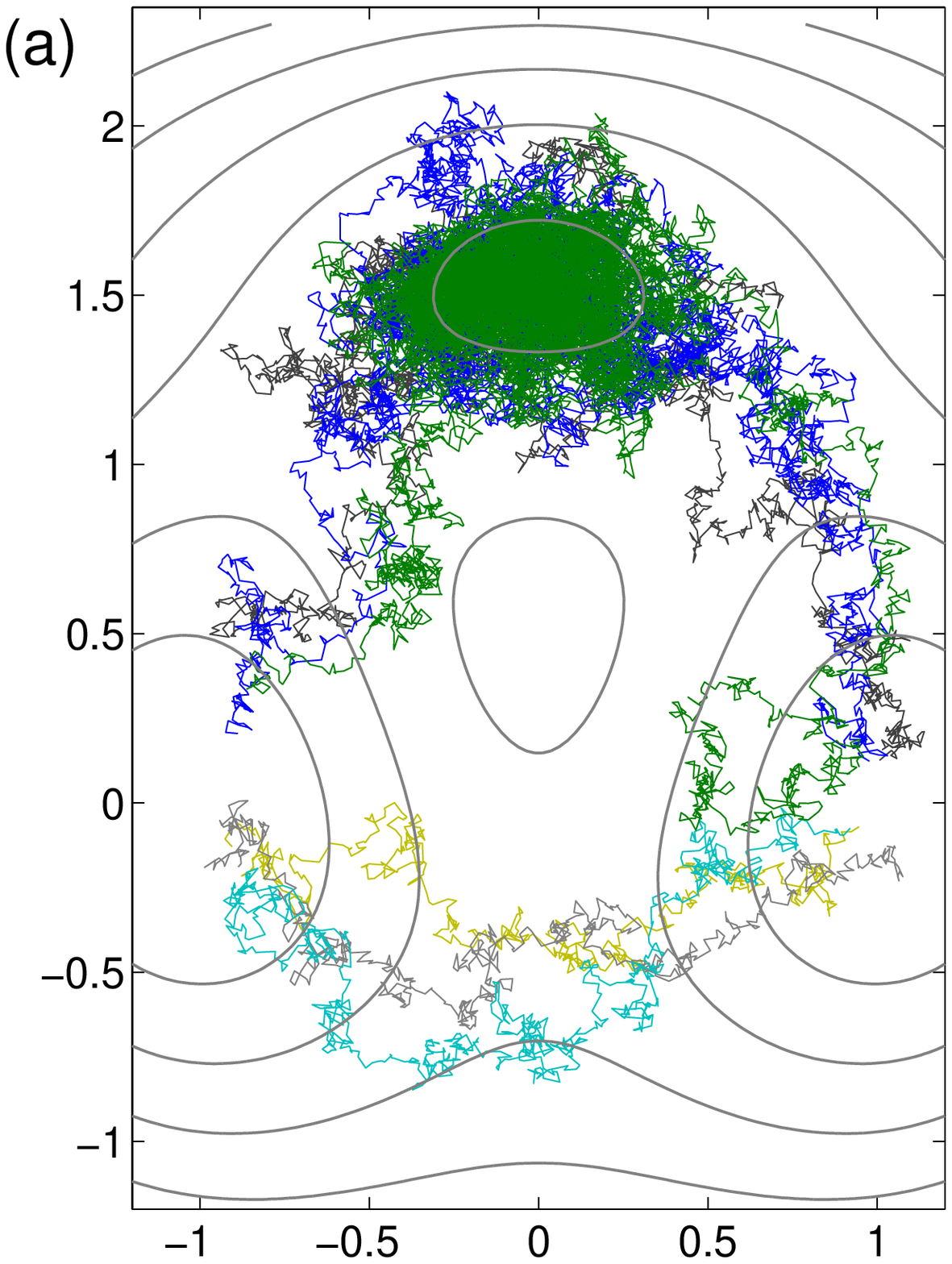}\includegraphics[height=5cm]{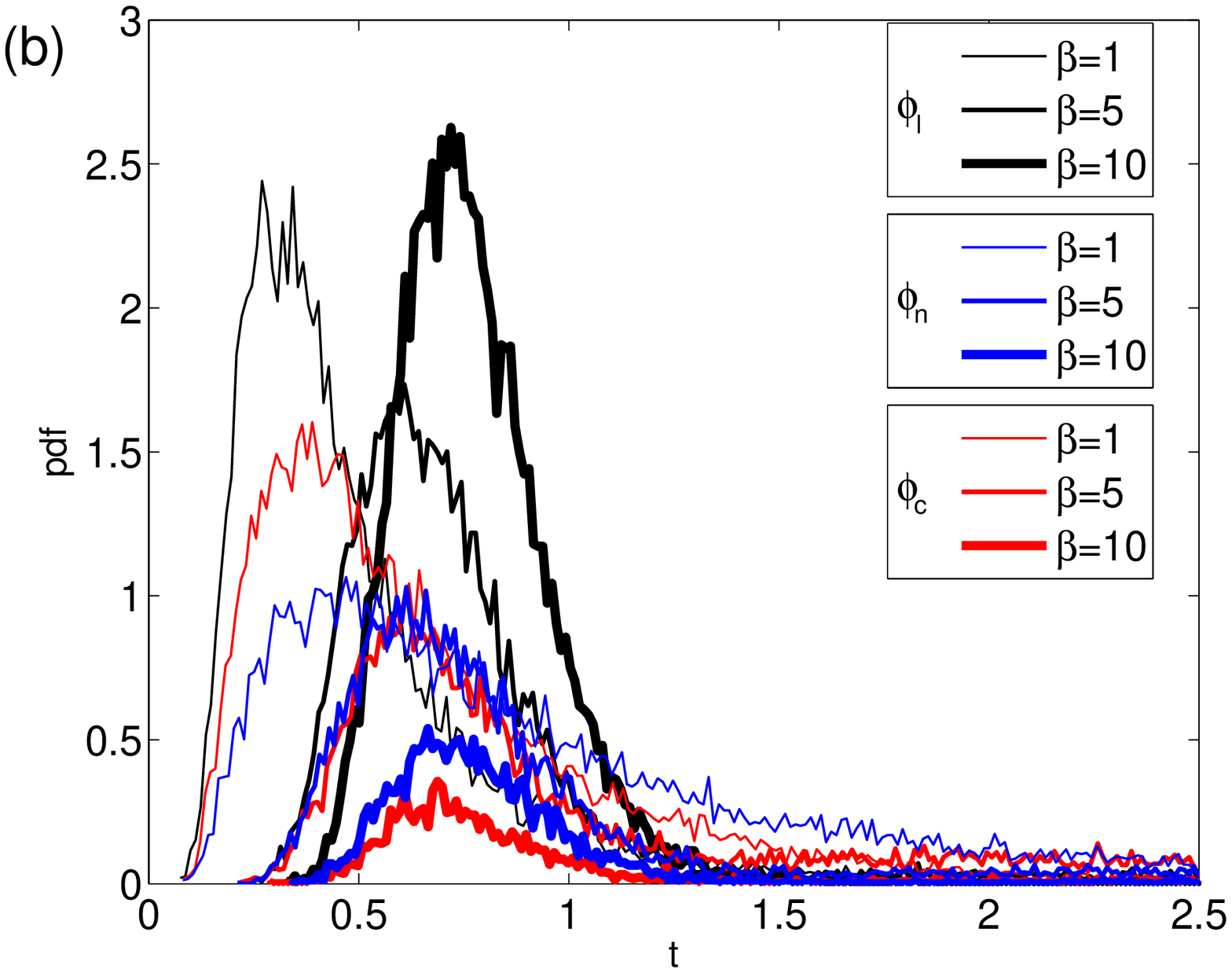}}
\caption{(a): Examples of trajectories corresponding to the limit between each decile of duration and
to the average duration, linear reaction coordinate at $\beta=5$ and $N=20000$. The
contour lines correspond to the  potential (Eq.~(\ref{tw})). (b) Example of
PDF of duration of reactive trajectories at $N=20000$ for increasing $\beta$, for all three reaction
coordinates. \emph{NB}: each
colour correspond to a reaction coordinate, each thickness of line to a temperature.}
\label{triple_}
\end{figure}

\begin{figure}
\centerline{\includegraphics[height=5.5cm]{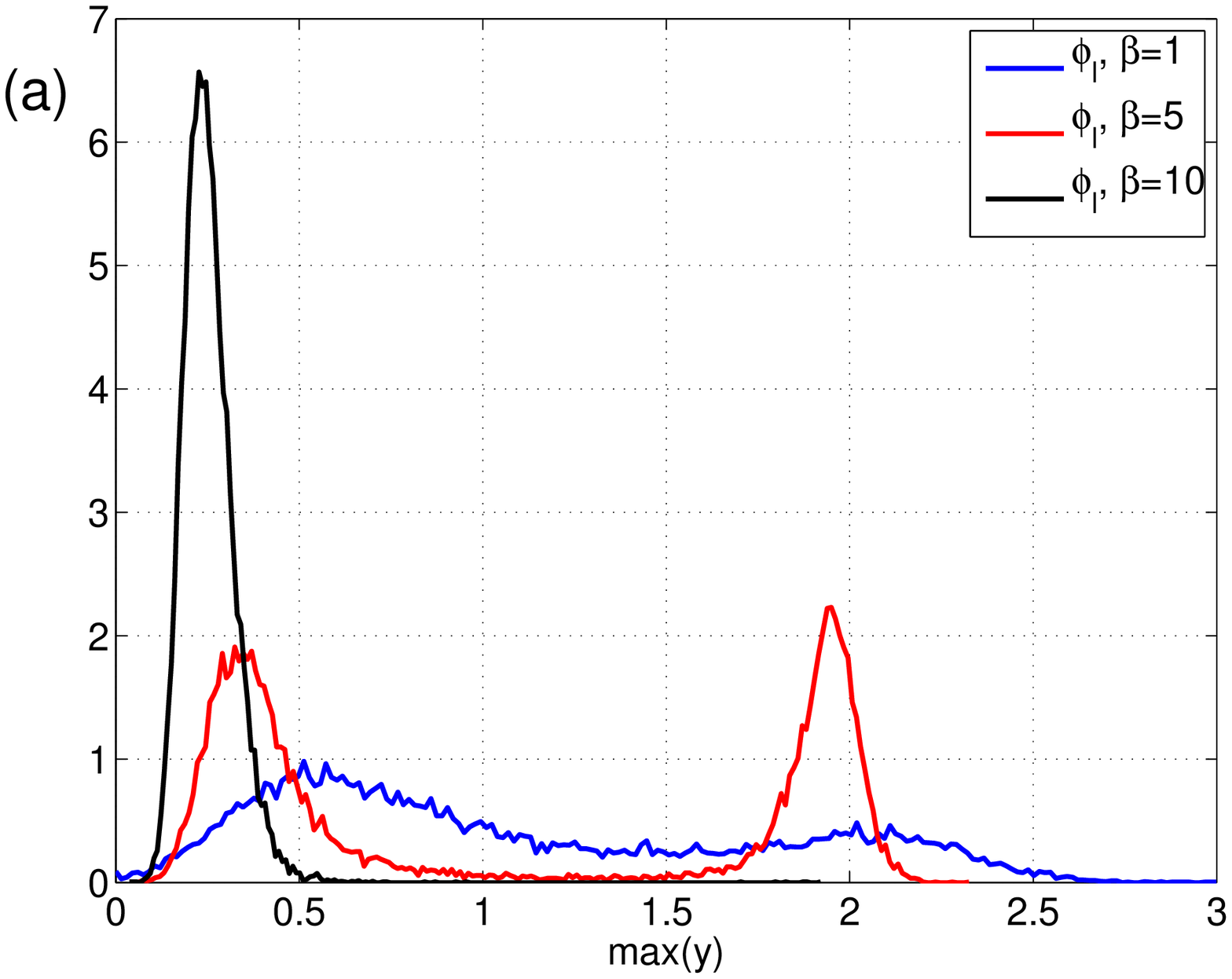}\includegraphics[height=5.5cm]{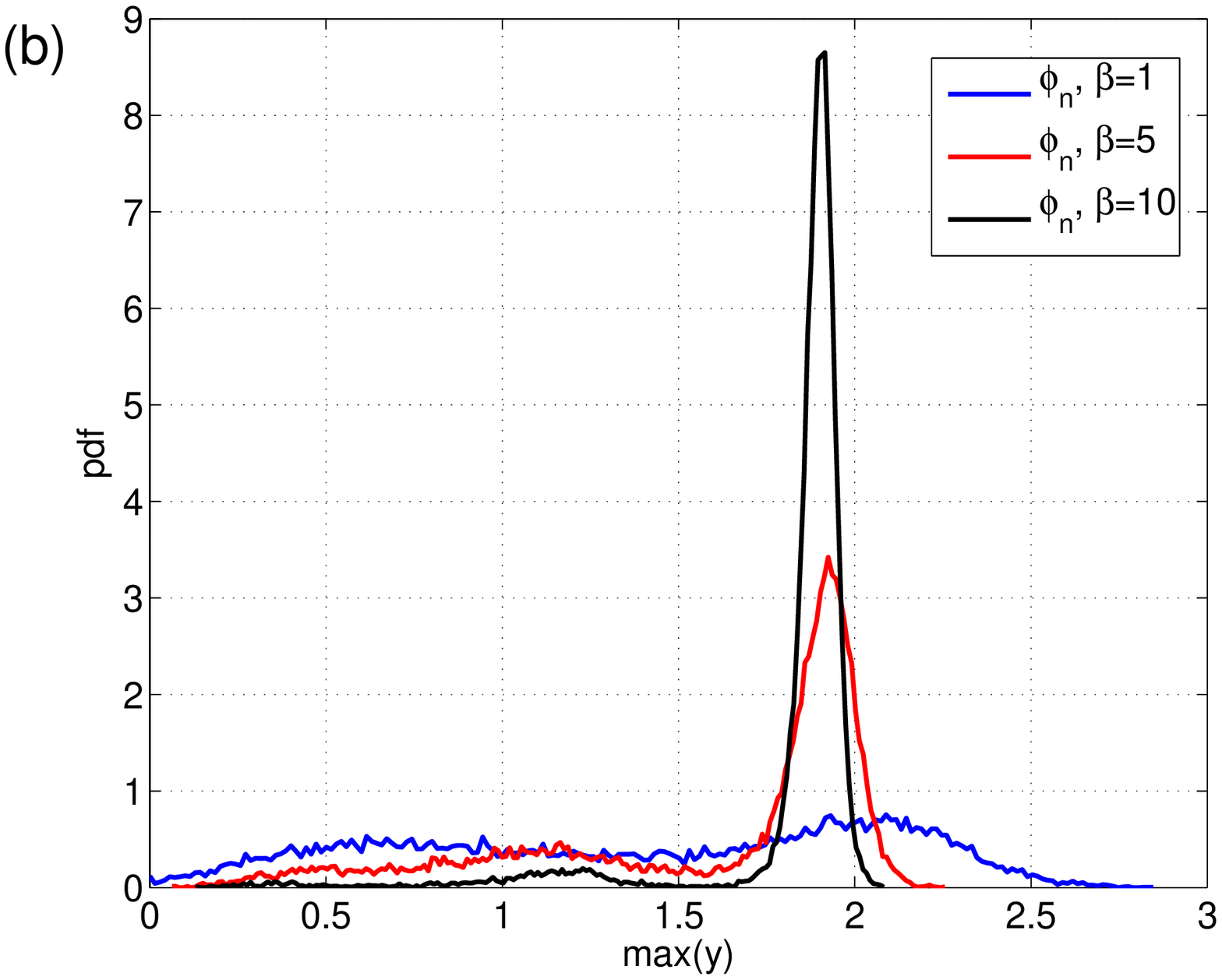}\includegraphics[height=5.5cm]{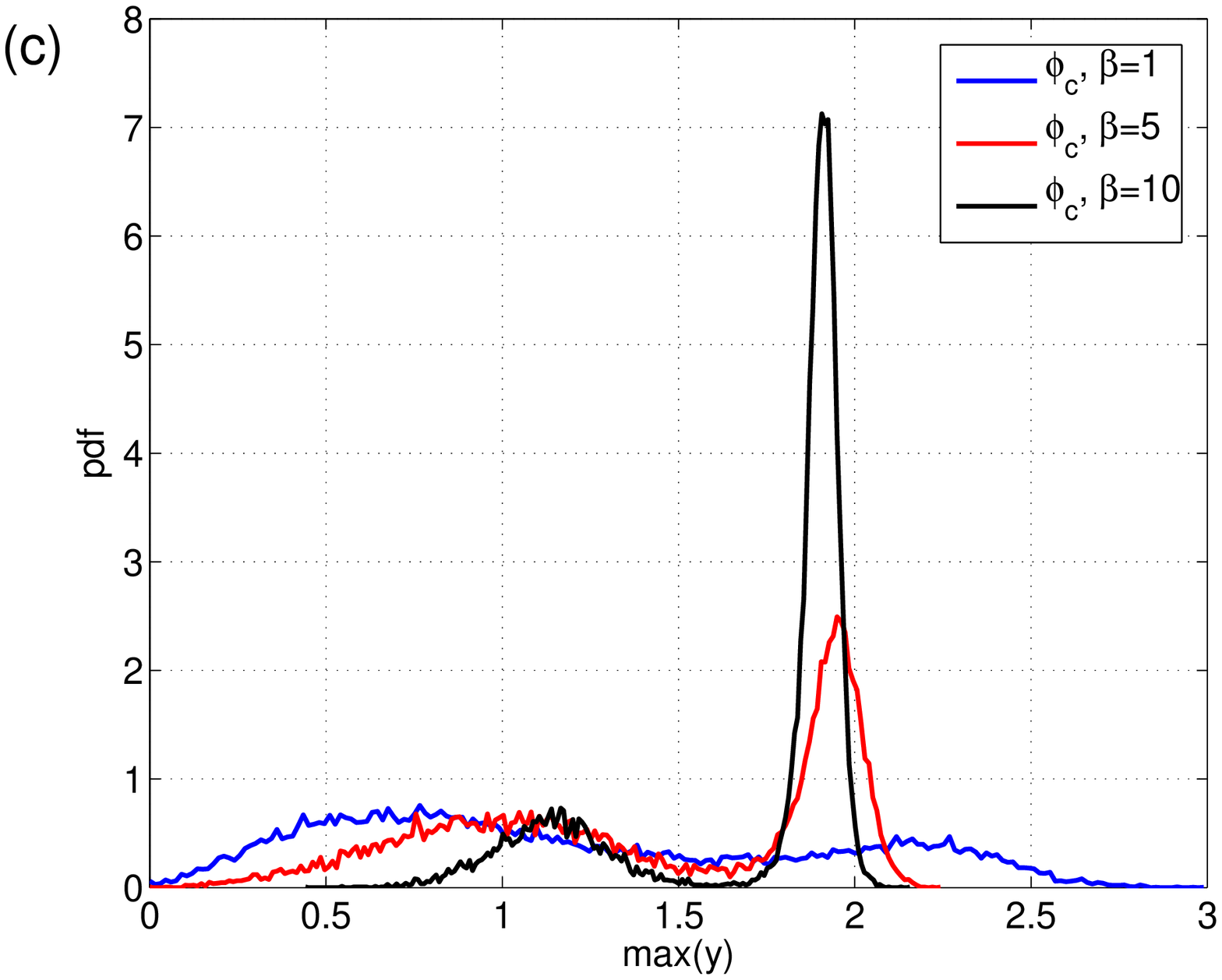}}
\caption{ Distribution of positions of the maximum value
of $y$ on each trajectory, for $\beta=1$, $5$, $10$, (a) for the linear reaction coordinate, (b) for the
norm reaction coordinate, (c) for the committor reaction coordinate.}
\label{triple}
\end{figure}

We focus here on the triple-well potential where a phase transition is observed at $\beta \simeq 5$.
We first examine the type of trajectories chosen for three different reaction coordinates.
For that matter, we choose $N=20000$, $dt=10^{-3}$, for $\beta=1$ (before the transition), $\beta=5$
(middle of the transition) and $\beta=10$ (after the transition). We show
in Fig.~\ref{triple_} (a) the two types of trajectories between the two minima:
the short duration ones using the lower channel, and the long duration ones
going through the metastable minimum \cite{cglp,VE}.
The effect of the choice of coordinates is striking when
examining the trajectories selected by the algorithm, distinguished by the maximum of $y$, $y_m$ on the
trajectory (Fig.~\ref{triple}), as well as the distribution of durations (Fig.~\ref{triple_} (b)).

The choice of the linear reaction coordinate tends to favor
short duration trajectories going to the lower part of
the domain. At large temperatures and near the phase transition, one observes both trajectories going to the
upper and lower part of the domain but as the temperature decreases, the upper channel
is not seen by the algorithm (see Fig.~\ref{triple} (a))
although it should in principle correspond to preferred trajectories. Note also that
the duration of reactive trajectories has short tails as seen in Fig.~\ref{triple_} (b).

The norm and committor coordinates have overall the same qualitative behavior (see Figs.~\ref{triple} (b),(c)).
One of the noticeable difference is that
the committor reaction coordinate does not suppress all trajectories going to the lower part even
in the small temperature regime. In both cases, the
distribution of durations has long tails (Fig.~\ref{triple_} (b)) which contrasts
with the linear coordinate short tail behavior. The result of the use of
the committor as a reaction coordinate is, as expected, in agrement with earlier computations \cite{cglp,VE}.
In fact, the committor seems to mix in an optimal way the two opposed tendencies
of the norm and linear coordinates and for all values of $\beta$.

Note that the examples with $\Phi_l$ presented here are typical. However, due to the large variance,
one can find rare, although possible, cases that do not differ much from
the two computations using $\Phi_n$ or $\Phi_c$ (not shown).

\subsection{Statistical behavior in the 2-D case}

The differences shown in the former section indicate a non-trivial behaviour of the algorithm. It can be seen
in the distribution of number of iterations $K$.
It is computed with repeated independent realisations (more than $10000$ in that case) of the algorithm for the same parameters
(here, $N=100$ and $\beta=5$). The Poisson distribution (\ref{poisson}) of the number of iterations,
in the Gaussian limit, has relative fluctuations and skewness:
\begin{equation}
\frac{\sigma}{m}=\frac{\sqrt{\langle K^2\rangle-\langle K\rangle^2}}{\langle K \rangle}=1\,,\,
S=\frac{\langle K^3\rangle-3\langle K\rangle\left(\langle K^2\rangle-\langle K\rangle^2 \right)
-\langle K\rangle^3}{\left(\langle K^2\rangle-\langle K\rangle^2\right)^{\frac{3}{2}}}=0\,.
\end{equation}
When examining the distribution of $K$ with different reaction coordinate (Fig.~\ref{dist_num} (c)),
one finds the following results (Tab.~\ref{tab1}).
\\
\begin{table}[!h]
\centerline{
\begin{tabular}{|r|c|c|}
\hline
$\Phi$ & $\sigma/m$ & S  \\
\hline
Committor & 1.6  & +0.05 \\
Linear    & 3.9  & -1.70 \\
Norm      & 10.7 & +2.30 \\ \hline
\end{tabular}}
\caption{Table summarising the values of the first cumulants of the PDFs of the number of iterations (triple well, $N=100$, $\beta=5$) for the three reaction coordinates.}
\label{tab1}
\end{table}
The distribution of $\hat{\alpha}$ is thus closer to the Poisson distribution (in the Gaussian limit)
when the committor is used. As we will see in the next section, the norm
coordinate improves greatly with respect to the linear reaction coordinate as $N$ is increased.

This result motivates a more detailed investigation of the statistical behavior of AMS on our 2-D models as a function
of the reaction coordinate and $\beta$ but also for large $N$. We try to keep
the number of independent realisations larger than 2000, although it
goes down to several hundreds (as used in \cite{cg07,pdm}) for the most expensive simulations
where $N=50000$.

\begin{figure}
\centerline{\includegraphics[width=6.5cm,clip]{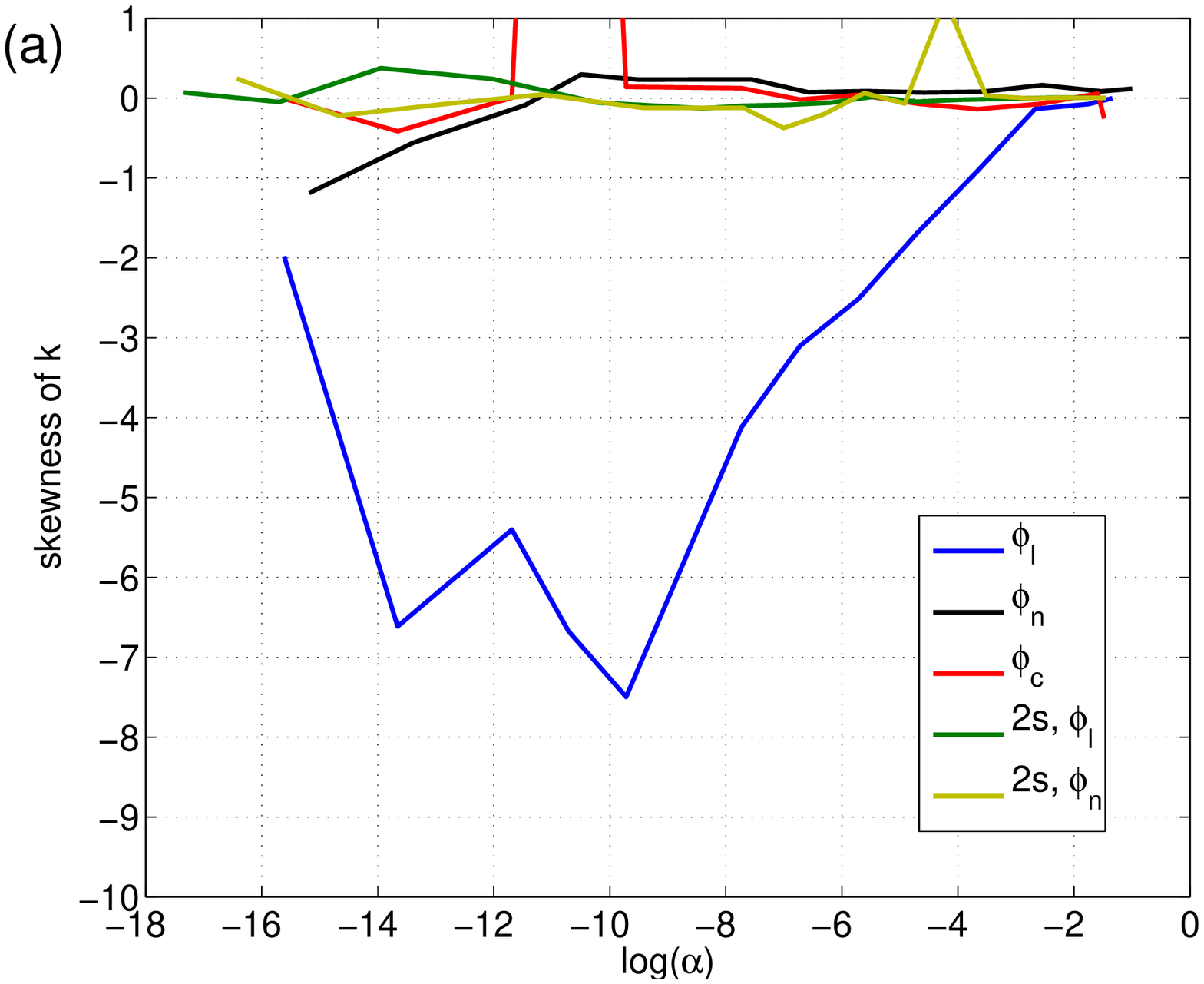}
\includegraphics[width=6.5cm,clip]{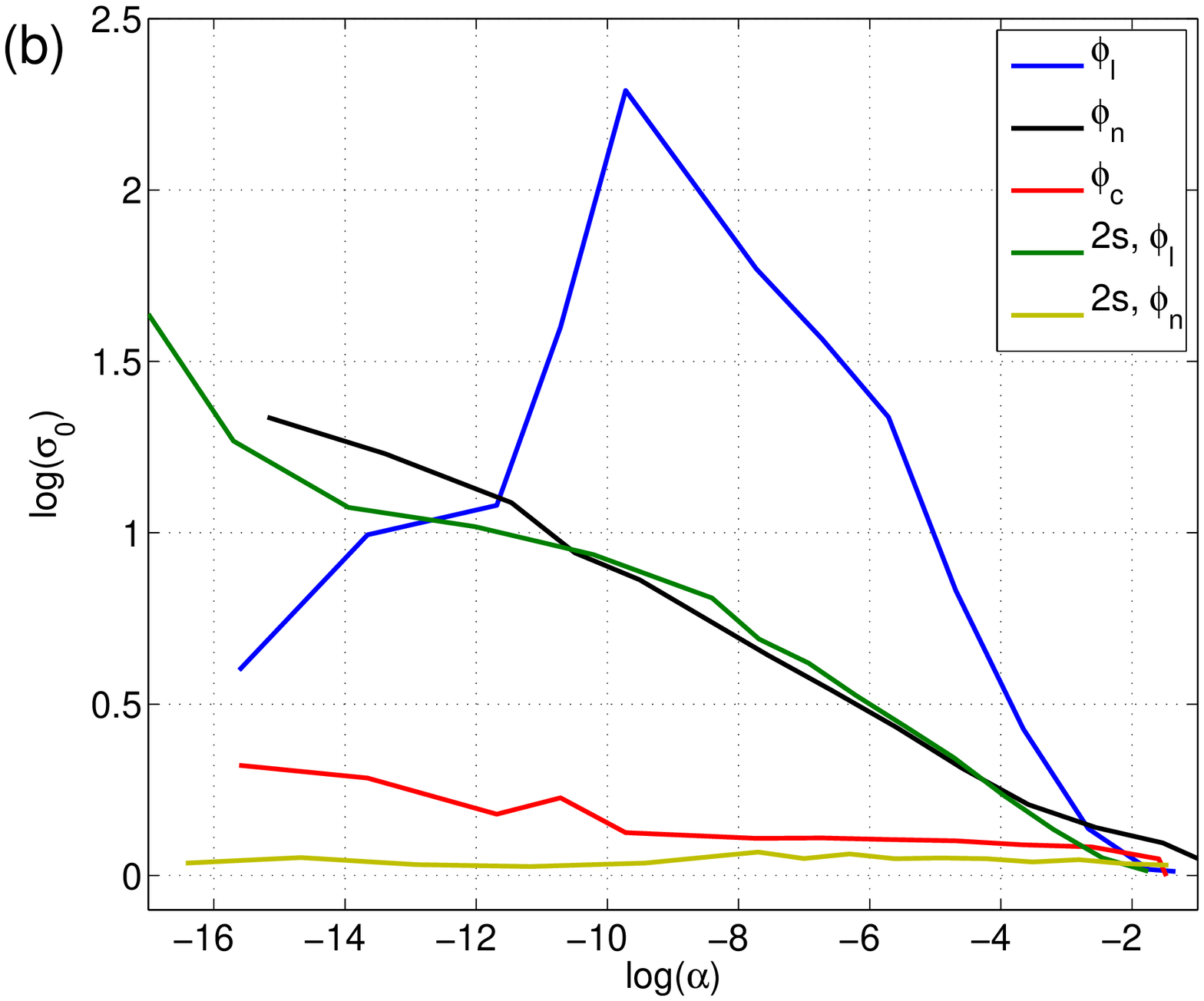}}
\centerline{\includegraphics[width=6.5cm]{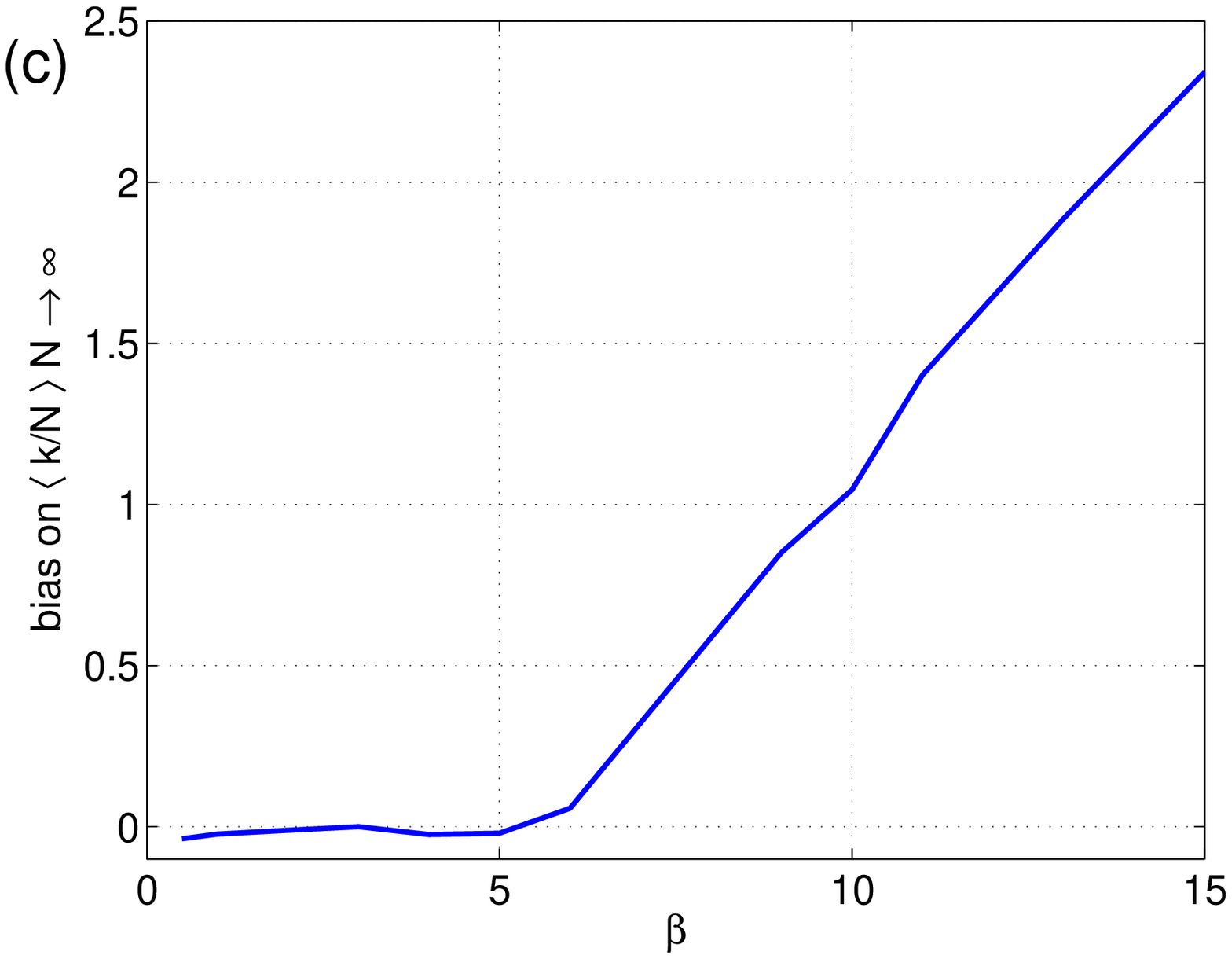}
\includegraphics[width=6.5cm]{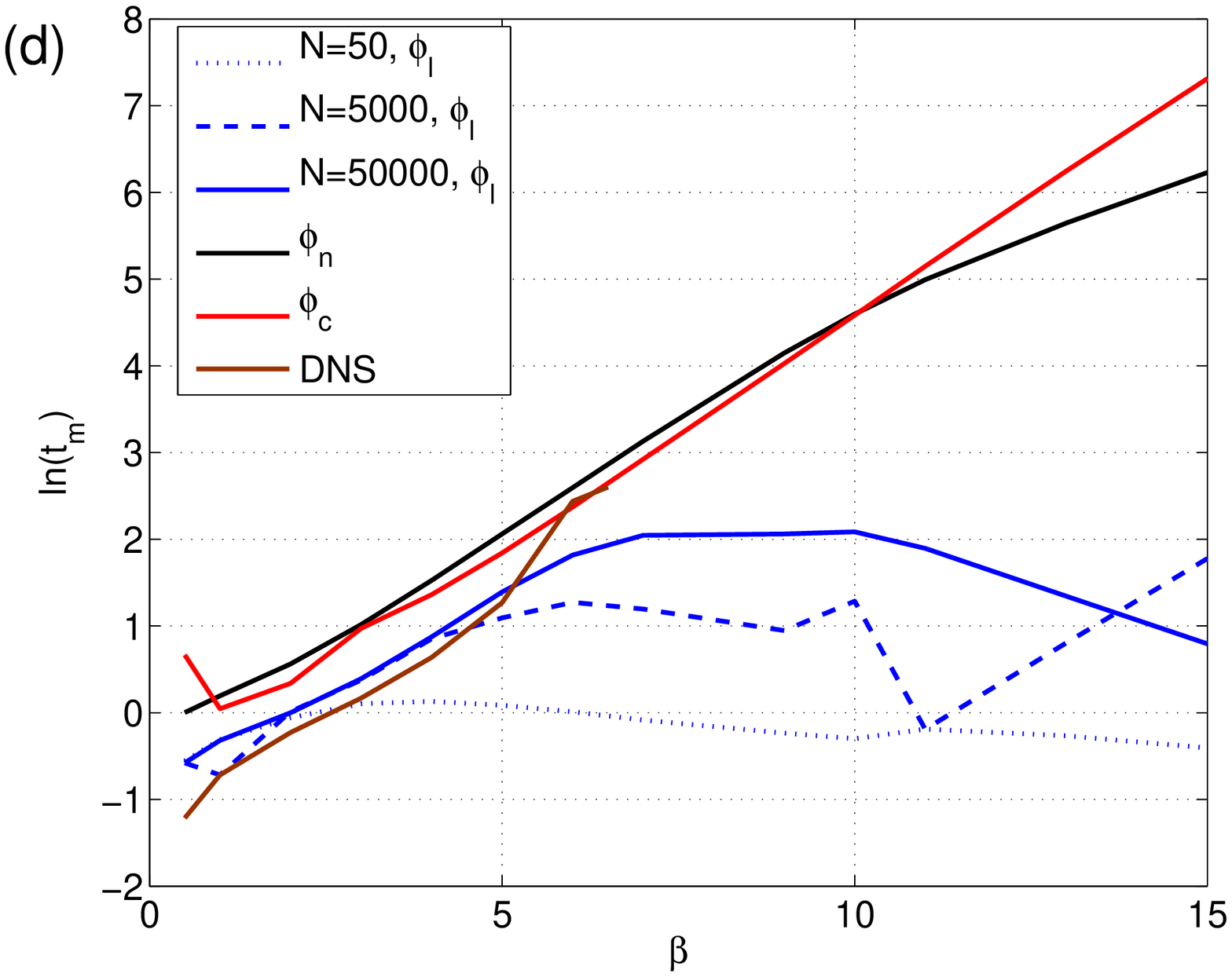}}
\caption{Comparison of statistics of the number of iterations of the algorithm between the triple-well potential
and the two saddles model through the transition as a function of $\ln(\alpha)$ for increasing $N$, and several reaction coordinates.  (a): Skewness of the crossing probability as a function of $\ln(\alpha)$ for several $N$, reaction coordinates.  (b): multiplied variance of the crossing probability as a function of $\ln(\alpha)$ for several $N$, reaction coordinates. (c) : bias on the average relative number of iterations $\langle k\rangle /N$ as a function of $\beta$.  (d): Average duration of reactive trajectories as a function of $\beta$ .}\label{trans_AMS}
\end{figure}

We now examine systematically the transition in this model using DNS and the algorithm with the three reaction
coordinates (Eq.~(\ref{coord_2c})). We run DNS for $\beta=0.5$ to $\beta=6.5$ and the algorithm for
$\beta=0.5$ to $\beta=15$. We use $dt=10^{-3}$. Computations using $dt=10^{-2}$ and $dt=10^{-4}$ (not shown here) gave the same results.
The number of clones goes from $N=10$, to $50000$ when the linear coordinate is chosen. We use only $N=2500$ otherwise.

We first compare the statistics of the number of iterations of the algorithm between the triple-well potential
and the two saddles model. For that matter, we compute a compensated variance:
$$
\sigma_0 \equiv -\frac{\sqrt{N}}{\alpha\sqrt{|\ln{\alpha}|}} {\rm Var}(K) .
$$
This quantity is equal to one for known theoretical predictions (see Eq.~(\ref{variance1})).
The relative skewness of the number of iterations in the perfect situation
should converge to zero in the asymptotic regime $N \to \infty$ (Fig.~\ref{trans_AMS} (a)).
On the other hand, we consider the convergence of the relative number of
iterations and of the estimate of the average durations of reactive trajectories.

The choice of the reaction coordinate has a striking effect on the skewness of the
distribution of number of iterations $K$ (blue curves in Fig.~\ref{trans_AMS} (a)) for the triple-well potential.
Choosing the linear reaction coordinate causes a negative asymmetry of the distribution
which is far greater than what
is observed with the norm or committor reaction coordinate. An
asymptotic regime has been reached: the skewness is independent of $N$ in the range
considered, even when $N$ is large enough to obtain good estimates
($\beta \lesssim 6$). Therefore in the case of the linear reaction coordinate, the data displayed here are averaged over $N$.
This asymmetry is negative as opposed to the positive asymmetry of the Poisson distribution.
Besides, it is much larger than that of the Poisson distribution in the
Gaussian limit. The skewness is minimum near the transition.
The use of the norm and committor reaction coordinate
leads to a positive skewness. Said skewness is small for all $N$ when using the committor reaction coordinate, while it decreases rapidly with $N$ when using the norm reaction coordinate. Indeed, it goes from the apparently large value of table~\ref{tab1} at $N=100$ to much smaller values if $N\gtrsim 250$. By contrast, the two saddles model keeps a very
small skewness of the distribution of number of iterations.

The reaction coordinate has also a clear effect on the compensated variance $\sigma_0$.
We display its decimal logarithm in figure~\ref{trans_AMS} (b) for a better appreciation of the orders of magnitude.
When using the linear coordinate, it has a behaviour comparable to the skewness. It is much larger than expected and is maximum at the phase transition. Again, this effect is independent of the number of clones
$N$: the data sampled using the linear reaction coordinate is averaged over $N$.
It eventually goes to $1$ at larger $\beta$. Meanwhile, using the norm coordinate leads to
no explosion, but a slow increase of $\sigma_0$ as the probability of crossing goes to zero.
Note that the rate of increase as a function of $\log(\alpha)$ is the same in the
triple-well potential using the norm coordinate $\Phi_n$ and in the two saddles model using the linear
coordinate $\Phi_l$. This indicates another generic behaviour of the algorithm: a slow
increase of the variance but a skewness remaining small.
This is more likely to be related to the loss of performance of the algorithm as $\beta$
is increased rather than to a poor rendering of the entropic switching during phase transition.
The committor reaction coordinate $\Phi_c$ yields
a compensated variance which is always close to one. This is also the case for the two saddles model
with the norm reaction coordinate.

We now look at the bias behavior as a function of the number of iterations for the triple-well potential only.
The error on the computation of $\alpha$ is examined \emph{via} a proxy, the average
relative number of iterations $\langle K\rangle/N$. It appears that this quantity gives an unbiased estimate
of $-\ln \alpha$ with variance $-\ln\alpha$ when the algorithm is either perfect or the committor is used
\cite{STCO}. In fact at an heuristic level, one has consistency with the central limit theorem (\ref{clt}) since
$ \tilde{\alpha} = (1-1/N)^K \simeq \exp(-K/N)$ for $N$ large.
Besides, the computation of this estimator yields more regular data than $(1-1/N)^K$.
As $N$ goes to infinity $\langle K\rangle/N$ converges
toward an asymptotic value which appears to have a finite bias relatively to
$\ln(\alpha)$ (computed by DNS for small $\beta$ or using the committor reaction coordinate for
large $\beta$). This bias is displayed in Fig.~\ref{trans_AMS} (c). One can see that
after a threshold value of $\beta\simeq 5$, in the middle of the transition range, it
grows nearly linearly. As a consequence, relative errors of more than $10$ are made in the estimation of $\alpha$, even for relatively large values of $\alpha$. Note that the slope is approximately $0.25$, which is very close to what one would have if $\Phi_l$ estimated the probability of having a transition through the high saddle $\propto \exp(-2.6\beta)$. One has $\alpha \propto \exp(-2.32 \beta)$, indeed, results like the Freidlin--Wentzell theory tell us that the crossing probability goes like $\exp(-\beta \Delta V)$ in potential systems \cite{ht}. Then $K/N$ estimates $2.6\beta$ instead of $2.32\beta$, when $\Phi_l$ is used. This would lead to a bias growing like $0.28\beta$. This shows that this is the effect of the selection
of trajectories going through the saddle by the linear coordinate.

The most striking effect of the choice of reaction coordinate is on the duration
of reactive trajectories (Fig.~\ref{trans_AMS} (d)). Both DNS and the use of the norm and committor
coordinates show an exponential increase of the average duration of the trajectory,
which corresponds to the mean first passage from the metastable minimum to the set $\mathcal{B}$. When using the linear reaction coordinate, the computed
trajectory durations behave quite differently as a function of $\beta$. The duration first follows the exponential growth up to a certain point,
then saturates and decreases toward an asymptotic behaviour which corresponds to trajectories
going through the saddle. The inverse temperature for which the systems leaves the expected
behaviour increases as $N$ increases.
The more clones are used, the more the algorithm produces long trajectories going through the metastable
minimum.  As a consequence, average durations computed using the linear reaction coordinate may very well
converge as $N\rightarrow \infty$. However, the rate of convergence decreases dramatically as
$\beta$ is increased. No matter how large is $N$, there exist a $\beta$ for which the algorithm will
eventually go back to selecting trajectories through the lower path.

\begin{figure}
\centerline{\includegraphics[width=6cm]{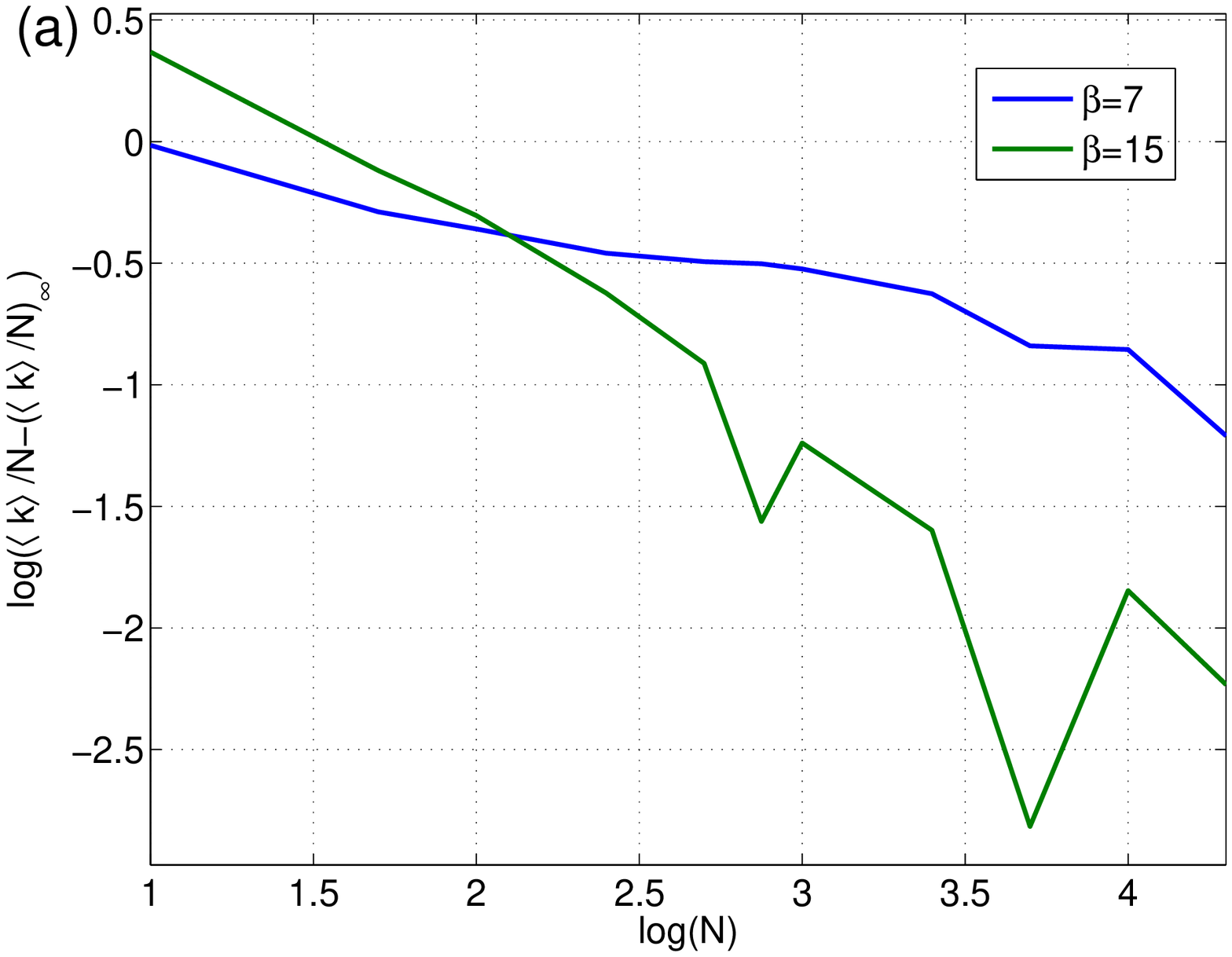}\includegraphics[width=6cm]{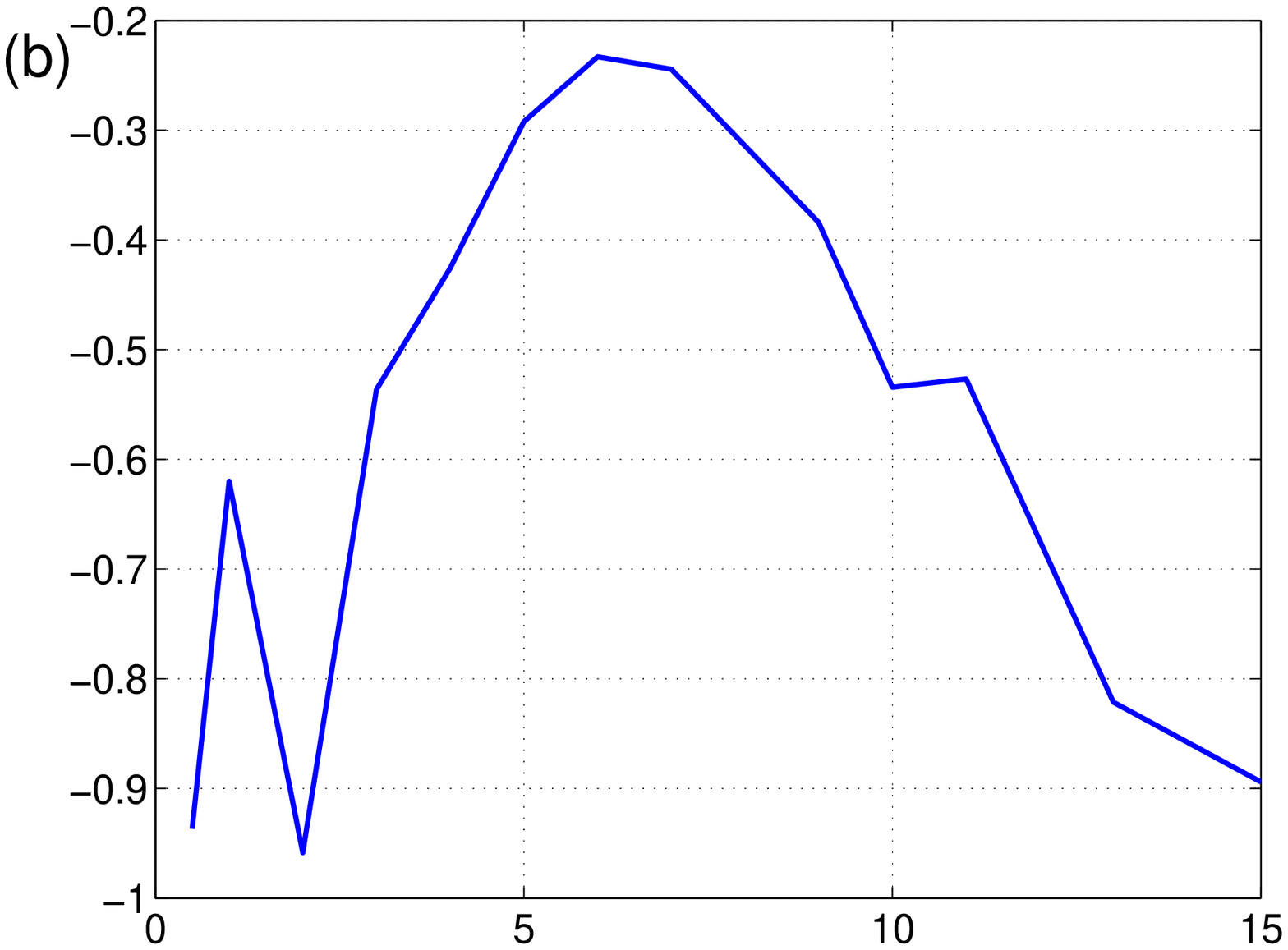}\includegraphics[width=6cm]{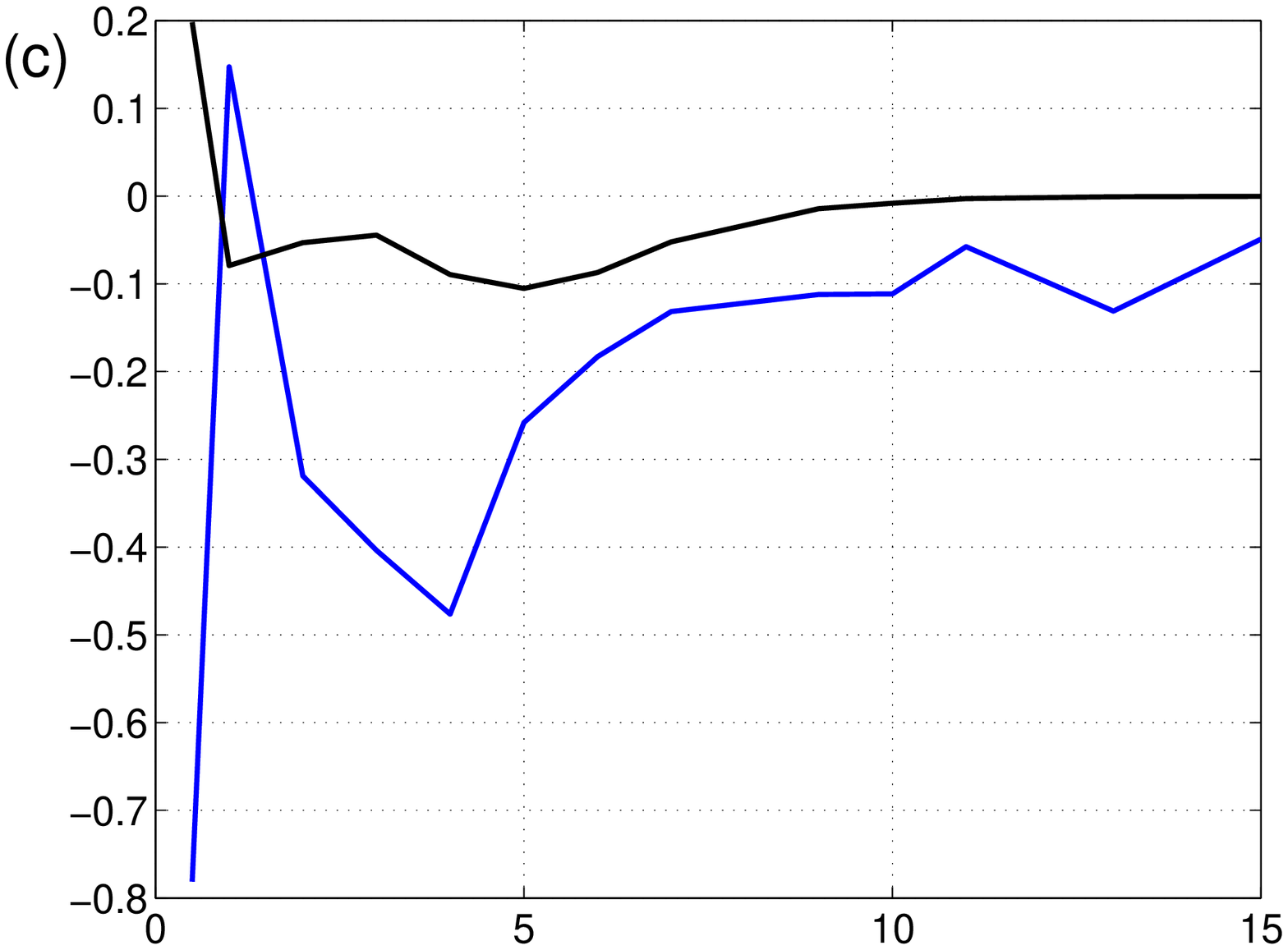}}
\caption{Convergence of the relative number of iterations ($\log(|\langle K\rangle/N-(\langle K\rangle/N)_\infty|)$) as a function of
$\log(N)$ for $\beta=7$ and $\beta=15$. Rate of convergence  as a
function of $\beta$ of (b):  the relative number of iterations $f_{\alpha}(\beta)$, (c): the average duration of reactive trajectories  $f_{\tau}(\beta)$ (blue : assumption of non sytematic bias, black : assumption of a systematic bias).}
\label{rates}
\end{figure}

This can be verified quantitatively by computing the rate $f(\beta)$ at which the $N$
dependent bias on the estimation of $\tau$ and $K/N$ converges toward $0$. This rate is defined by :
\begin{equation}
\langle\hat{\tau}_N\rangle-\varpi\propto N^{f_\tau(\beta)}\,,\,\left\langle\frac{K}{N}\right\rangle-\lim_{N\rightarrow \infty}\left\langle\frac{K}{N}\right\rangle\propto N^{f_\alpha(\beta)}\,,
\end{equation}
with $\varpi=\tau$ or $\hat{\tau}_{\infty}$ depending on how the exact value is estimated:
$\tau$ is computed using the committor reaction coordinate when DNS is not available and
$\hat{\tau}_{\infty}$ is the asymptotic limit $N\rightarrow \infty$ approximated by $N=50000$.
For the average duration of reactive trajectory, we compute the difference with both
the asymptotic value and the exact value (the biases are respectively termed relative and absolute) because it is
not clear whether $\hat{\tau}$ has a systematic bias or not. There are two distinct behaviours for the rate of
convergence of $\tau$ and $K/N$. The rate of convergence of the relative number of iterations $K/N$
is always finite (Fig.~\ref{rates} (a)), and we can distinguish two regimes : one where it is small,
typically in the middle of the transition range, and one where it is close to $-1$, before and after
the transition range (Fig.~\ref{rates} (b)). Similarly, $f_{\tau}$ departs from the rather
quick convergence regime exhibited numerically in section~\ref{se_norm} (Fig.~\ref{stat_temps} (b)).
However, it goes rapidly to zero, whether we assume $\tau$ has a systematic bias or not. This seems to
indicate that in the case of $\tau$, the convergence is always logarithmic (Fig.~\ref{rates} (c)).  Note that the use of $dt=10^{-2}$
shows that this effect is in fact independent of $dt$: the convergence rates (not indicated in the curves) are approximately equal to those at $dt=10^{-3}$.

In all cases, no improvement is seen when using the linear reaction coordinate and decreasing $dt$ to $10^{-4}$, for $N=2500$.
The convergence toward paths going through the metastable minimum occurs only when $N$ is increased.

\subsection{A three-level model}

We aim at describing the effect of the selection of reactive trajectories
by the linear reaction coordinate on the estimate of the average duration.
For that matter, we propose a simplified model with three states which mimics the behavior of
the triple-well potential system, and solve it analytically. We will aim at extracting a convergence
rate of the algorithm in that case.

\subsubsection{Definition}

We decompose the phase space in three states $1$, $2$, $3$. The states $1$ and $2$ correspond
to the two global minima of the triple-well potential $\mathcal{A}$ and $\mathcal{B}$, and the state $3$ corresponds
to the metastable minimum $\mathcal{D}$. Let $\rho_{1,2,3}$
be the probability of presence in those states, and $R=(\rho_1,\rho_2,\rho_3)$. We have $\sum_iR_i=1$,
the normalisation of probabilities. The dynamics of $R$ around equilibrium is given by a transition matrix \cite{VK}:
\begin{equation}
\frac{dR}{dt}=MR \,,\,M=
\left( \begin{matrix} -P_{1\rightarrow 2}-P_{1\rightarrow 3}& P_{2\rightarrow 1}&P_{3\rightarrow 1}
\\ P_{1\rightarrow 2}&-P_{2\rightarrow 1}-P_{2\rightarrow 3}& P_{3\rightarrow 2}
\\ P_{1\rightarrow 3}& P_{2\rightarrow 3}&-P_{3\rightarrow 2}-P_{3\rightarrow 1} \end{matrix}\right)\,.
\end{equation}
The $P_{i\rightarrow j}$ are the rate of transition between states $i$ and $j$.
One obviously has $\sum_{i}M_{ij}=0$, so that $d(\sum_{i}R_i)/dt=0$: the normalisation of probabilities is conserved.
Due to detailed balance, one has
 \begin{equation}
 P_{1\rightarrow 3}=P_{2\rightarrow 3}\equiv B\,, \,P_{3\rightarrow 1}=P_{3\rightarrow 2}\equiv C\,,\, P_{1\rightarrow 2}=P_{2\rightarrow 1}\equiv A\,.
 \end{equation}
The time evolution is obtained by diagonalising $M$ : $M=Q^{-1}DQ$, and using the initial condition, \emph{i.e.}:
\begin{equation} \frac{dR}{dt}=Q^{-1}DQR\,,\, Y=QR\,,\, R=Q^{-1}Y\Leftrightarrow \frac{dY}{dt}=DY\,.\end{equation}

We can adapt this framework to describe reactive trajectories. We set $2$ as an absorbing state, \emph{i.e.}
$P_{2 \to 3} = P_{2 \to 1} = 0$. Coming back to $\mathcal{A}$ is impossible, and transition from $\mathcal{D}$ to $\mathcal{B}$ going through the saddle between $\mathcal{D}$ and $\mathcal{A}$ and then through the saddle between $\mathcal{A}$ and $\mathcal{B}$ (avoiding $\mathcal{A}$) are extremely unlikely, so we set $P_{3\rightarrow 1}=0$. The Markov Matrix reads:
\begin{equation}
M=\left(\begin{matrix}-(A+B)&0&0 \\A &0& C \\ B&0 &-C \end{matrix} \right) \,.
\end{equation}
Note that, in that process, detailed balance is broken: this is because we consider only reactive trajectories. This modification is very similar to applying the h-transform of Doob in a Langevin equation \cite{doob}. This type of modification is in the same spirit as considering the distribution around a fluctuation trajectory, which is a solution of the Euler--Lagrange problem derived from the Freidlin--Wentzell theory \cite{ht}.
We use the inverse of the reactive trajectory duration for $A$, $B$ and $C$. The rate $A$ and $B$ correspond to reactive
trajectories out of the set $\mathcal{A}$. We can take either $A\propto B\propto 1/\beta$ or $1/\ln(\beta)$
(see \cite{tps}).
The former is more tractable analytically, and the later is more realistic. However, both will yield the same
qualitative behavior.
They have to decrease slowly, relatively to $C\propto e^{-\beta}$. This rate of transition corresponds to the mean first passage out of the metastable minimum \cite{ha,VK}. Indeed,
trajectories from $\mathcal{D}$ to $\mathcal{B}$ correspond to first passages which are not reactive,
the system is allowed
to fail to cross the saddle and to go back to $\mathcal{D}$ until it reaches $\mathcal{B}$.
We choose to balance the rate of probability $A$ and $B$, as a consequence, trajectories go indiscriminately
through both paths. We will see that it has no effects on the average durations.

\subsubsection{Results}

\begin{figure}
\centerline{\includegraphics[width=7cm]{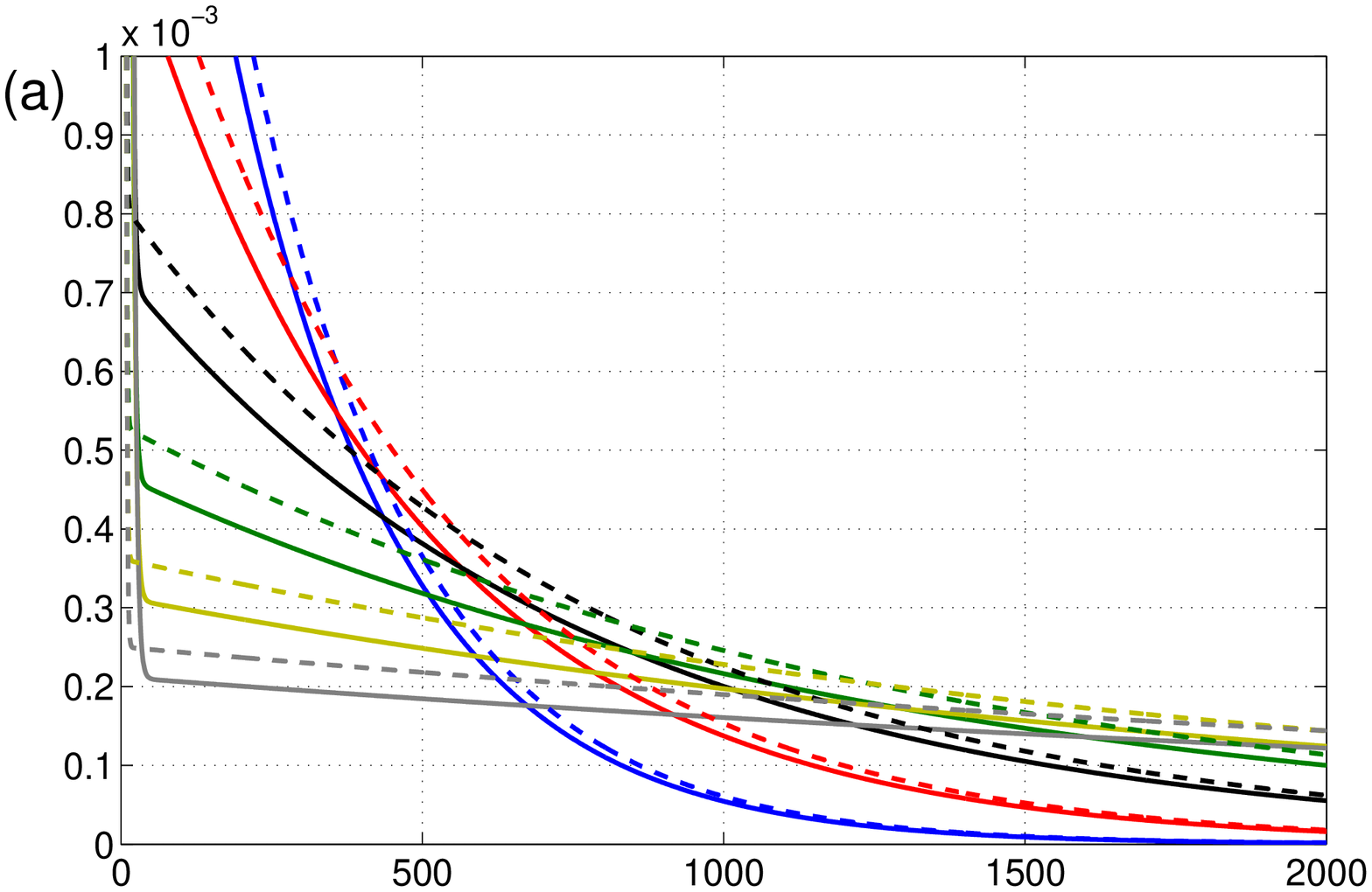}
\includegraphics[width=7cm]{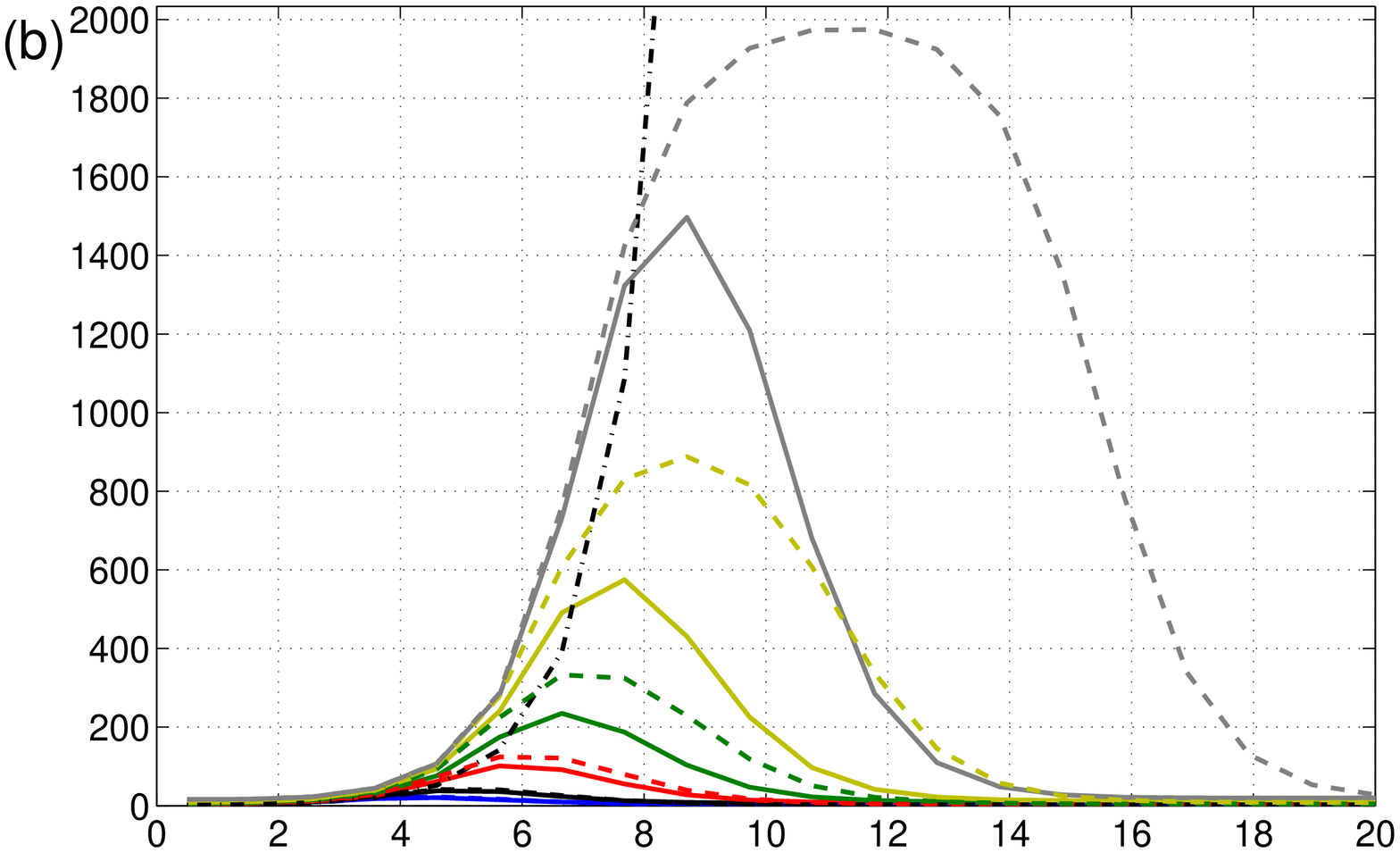}}
\centerline{\includegraphics[width=7cm]{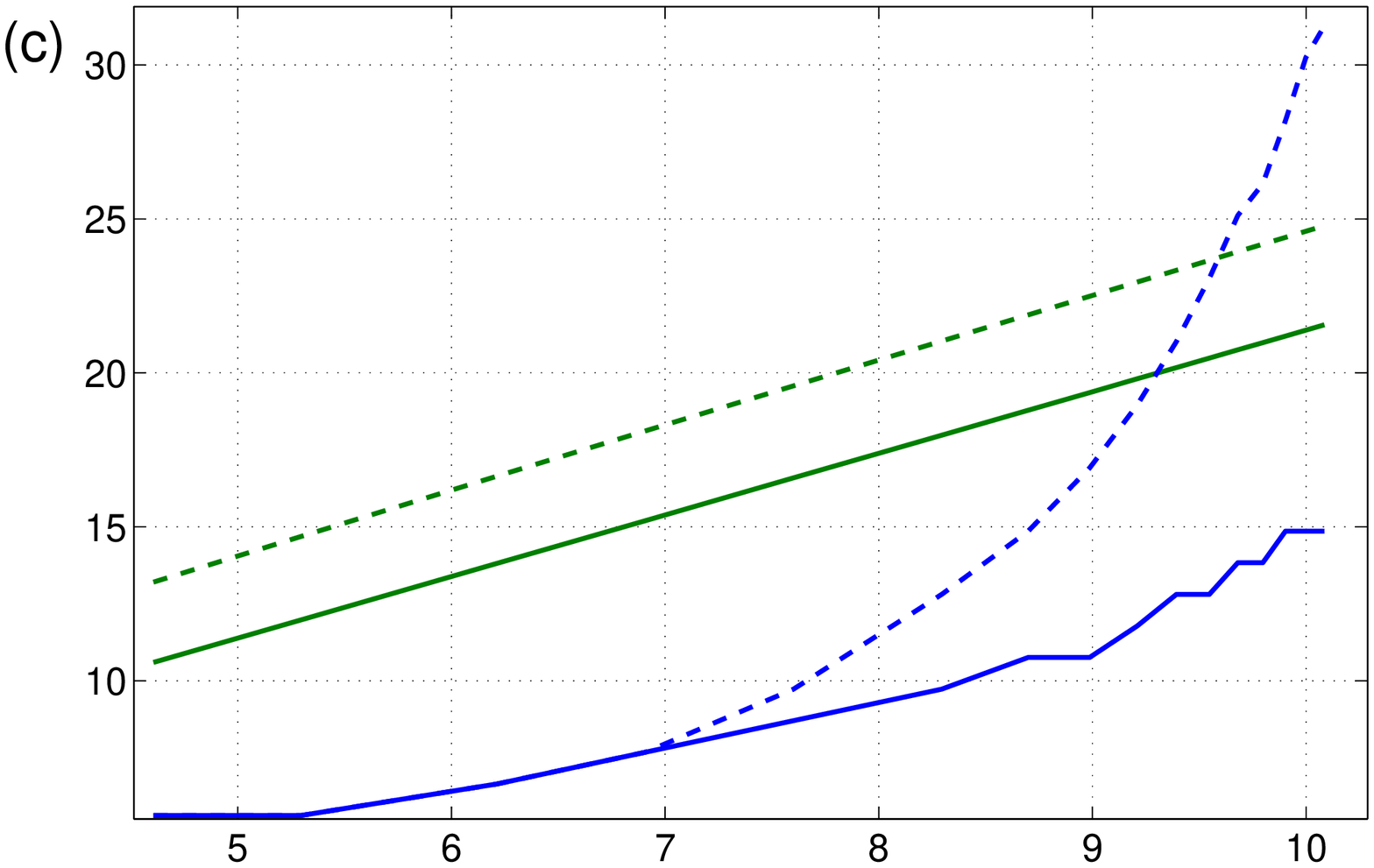}
\includegraphics[width=7cm]{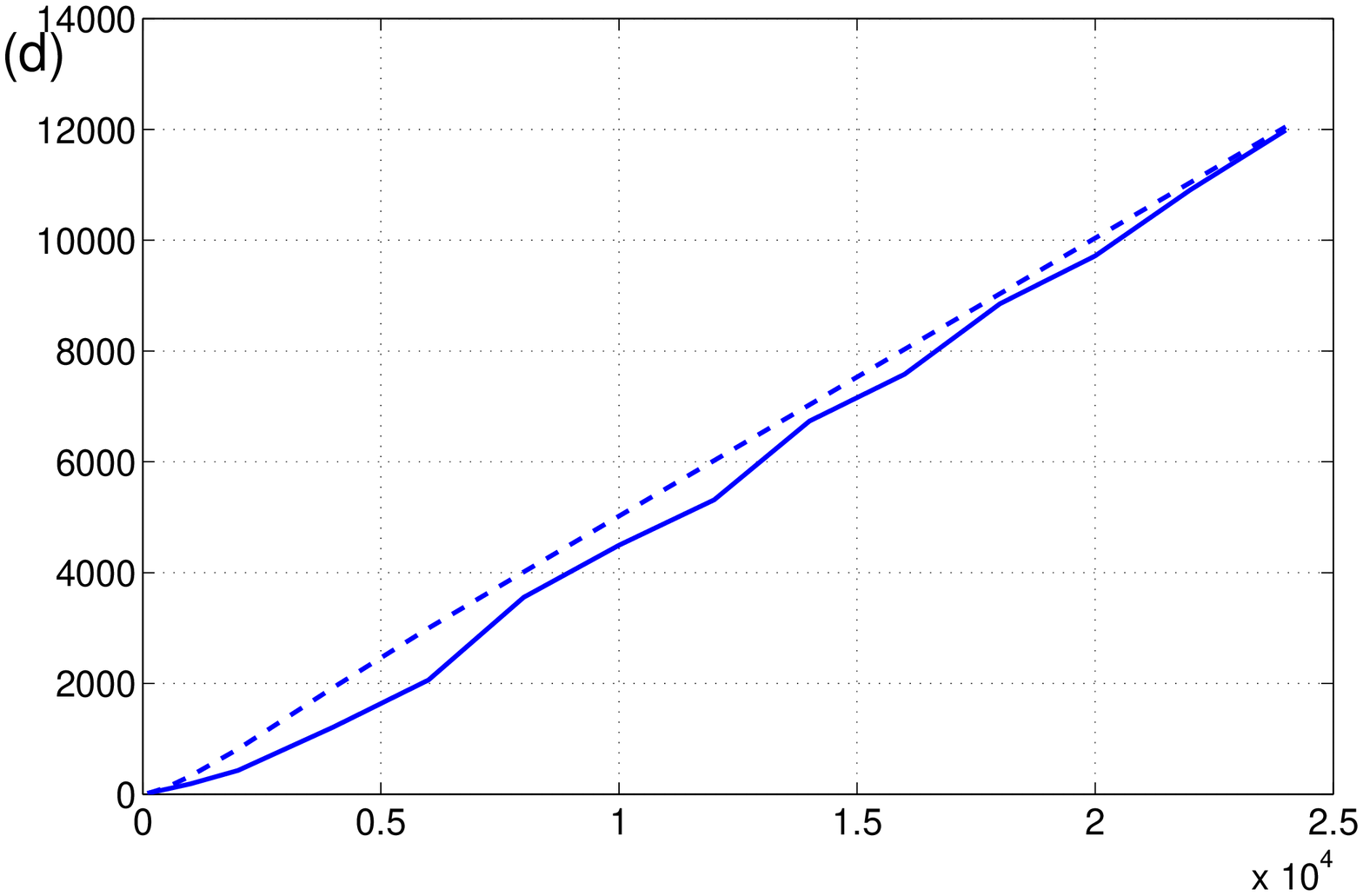}}
\caption{(a): PDF of the durations of reactive trajectories for several $6.1\le\beta\le8.6$ ($\Delta \beta=0.5$ increments) around the change of
regime. (b): average transition duration as a function of $\beta$ for several values
of the cut-off $\Lambda=100$ (blue), $200$ (black), $500$ (red), $1000$ (green), $2000$ (yellow), $4000$ (gray). The black dot-dash line is the analytical value. (c): value of $\beta$ at the maximum (blue) and analytical value of the inflexion
of $\langle \tau\rangle(\beta)$ (green) as a function of $\ln(\Lambda)$. (d): value of the maximum over $\beta$ of
the average duration as a function of $\Lambda$. Full line $A=B=1/\beta$, dashed line $A=B=1/\ln(\beta)$}
\label{fig}
\end{figure}

We compute the eigenvalues of $M$ defined by
\begin{equation}\det(M-\lambda Id)=-\lambda(C+\lambda)(A+B+\lambda)\Rightarrow
\lambda_1=0\,,\, \lambda_2=-(A+B)\,,\, \lambda_3=-C\,.
\end{equation}
They correspond respectively to the absorbing state $2$, the exit out of $1$ and the exit out of $3$.
The matrices $Q$ are:
\begin{equation}
Q^{-1}= \left( \begin{matrix} 0&1&0 \\[2mm] 1&\frac{BC+A(C-(A+B))}{(A+B)(A+B-C)}&-1\\[2mm] 0& \frac{B}{(C-(A+B))}&1 \end{matrix} \right) \,,\,  Q=\left(\begin{matrix}1 &1 &1 \\[2mm] 1 & 0 & 0 \\[2mm] \frac{B}{A+B-C} & 0& 1\end{matrix}\right) \,.
\end{equation}
We then find the matrix of transition, which gives $R(t)=T R(0)$ for an initial condition $R(0)$:
\begin{equation}
 g_1=\frac{B}{A+B-C}\,,\, g_2=\frac{BC+A(C-(A+B))}{(A+B)(A+B-C)} \Rightarrow   T(t)=\left(\begin{matrix}e^{-(A+B)t} & 0&0 \\[2mm]1-g_1e^{-Ct}+g_2e^{-(A+B)t} &0 &\left(1-e^{-C t} \right) \\[2mm]g_1\left( e^{-Ct}-e^{-(A+B)t}\right) & 0&e^{-C t}\end{matrix} \right) \,.
\end{equation}
The probability of transition in a time $t$ is the distribution function,
or cumulated density in physics literature, of the transition duration.
Indeed the probability of being in $\mathcal{B}$ at $t$ is the probability of having a transition duration
$\tilde{t}\le t$.
One then finds the PDF of transition durations from $1$ to $2$ by taking the derivative of $T_{2,1}$, which yields
:
\begin{equation}
 T_{12}'(t)=\frac{BC}{A+B-C}\left(e^{-Ct}-e^{(A+B)t} \right)+Ae^{-(A+B)t}=\frac{BC}{A+B-C}e^{-Ct}+\frac{(A+B)(A-C)}{A+B-C}e^{-(A+B)t}\,.
\end{equation}
It contains the transition through both passages.
From this, we extract the mean duration of a reactive trajectory $\tau = \int_0^\infty t T_{12}'(t) dt$:
\begin{equation}\tau= \frac{1}{A+B}\left(1+\frac{B}{C} \right)\,.\end{equation}
We display examples of $T_{12}'$ for several values of $\beta$ in Fig.~\ref{fig} (a).
The simplified dynamics does not contain a minimum duration for the reactive trajectory,
however, it captures the qualitative change of shape of the PDFs very well, when
the duration $1/C$ becomes large (see Fig.~\ref{triple} (c,d) for a comparison).
The exponential growth of $\tau$ can be seen in Fig.~\ref{fig} (b).

We want to investigate the effect of the number of clones when using the linear
coordinate $\Phi_l$. PDF of trajectories durations (Fig.~\ref{triple} (c,d)) and processed results of
the AMS (Fig.~\ref{trans_AMS} (d)) show that the proportion of long trajectories among those computed is
an increasing function of $N$. That is to say that the distributions of reactive trajectories have longer and longer tails as $N$ is increased. Qualitatively, it means that for low $N$ reactive trajectories tend
to favor the lower channel ${\cal A}$ to ${\cal B}$ with short duration. Therefore, we model the effect of
$N$ when $\Phi_l$ is used through a cut-off on the duration of trajectories $\Lambda$ as defined below.
It is a growing function of $N$. In practice the effect of $N$ is more continuous, but this constitutes a
good first approximation.
\begin{equation}
\langle f\rangle_{\Lambda}=\int_{0}^\Lambda f(t)\rho(t) {\rm d}t\,,\, \tau_\lambda=\int_{0}^\Lambda tT_{12}'(t) {\rm d}t\,, \,  \int_0^{\Lambda}te^{-at}{\rm d}t=\frac{1}{a^2}\left(1-(1+a\Lambda)e^{-a\Lambda}\right)\,.
\end{equation}

We compute the average including the cut-off. We display the results for both $A=B=1/\beta$
and $A=B=1/\ln(\beta)$ in figure~\ref{fig} (b).
The introduction of a cut-off captures very well the difference of behaviour between
the DNS and the AMS using the norm and committor coordinates, and the AMS with the linear
reaction coordinate (Fig.~\ref{trans_AMS} (d)).
As the cut-off is increased, the average duration $\tau$ in the model sticks to the actual behaviour for larger and larger $\beta$, then saturates and goes
back to the duration of transition through the saddle. Indeed, the cut-off on durations suppresses the trajectories
going through the metastable minimum.

For an analytical treatment, one
can make some simplifications, since $(A+B)\Lambda$ is always large. The average duration then reads:
\begin{equation}\tau_\Lambda= \frac{(A-C)C+B(A+B)\left(1-(1+C\Lambda)e^{-C\Lambda} \right)}{C(A+B)(A+B-C)}\,.\end{equation}
We expand it for $C$ small
\begin{equation}\tau_\Lambda=\frac{(A-C)C+B(A+B)\left(\frac{C^2\Lambda^2}{2}+O(C^3\Lambda^3)\right)}{C(A+B)(A+B-C)}=\frac{1}{A+B}\left( \underbrace{\frac{A}{A+B}}_{= cte}+CB\left(\Lambda^2-\frac{1}{(A+B)^2}\right)\right)+O(C^2)\,.\end{equation}

For $\beta$ large enough, and $C$ small, the average of the transition duration can be simplified:
\begin{equation}
\tau_\Lambda^{{ \rm lin}}=\frac{\beta}{2}\left(\frac{1}{2}+\frac{1}{\beta e^{\beta}}\Lambda^2 \right) \,,\, \tau_\Lambda^{{ \rm log}}=\frac{\ln(\beta)}{2}\left(\frac{1}{2}+\frac{1}{\ln(\beta) e^{\beta}}\Lambda^2 \right)
\end{equation}
which is minimum (\emph{i.e.} goes back to a growth like $\beta$) at:
\begin{equation}\beta_i^{\rm lin}=2\ln(\Lambda)+\ln(2)\,,\end{equation}
which increases with the cut-off.  If the growth of the duration of the reactive trajectories is
very slow, $A=B=1/\ln(\beta)$, the position of the minimum $\beta_i^{\rm log}$ verifies :
\begin{equation}
\beta_i^{\rm log}-\ln(\beta_i^{\rm log})=2\ln(\Lambda)+\ln(2)\Rightarrow \beta_i^{\rm log}=\beta_i^{\rm lin}+\ln(\beta_i^{\rm lin})/(1-1/\beta_i^{\rm lin})\,, \label{log}
\end{equation}
approximated at the first order in the difference to $\beta_i^{\rm lin}$.
This yields a simple upper bound on the inverse temperature $\beta$ for which
the model leaves the expected behaviour of the exponential growth. We can compare the computed positions
of inflexions to the maximum of $\tau_\Lambda(\beta)$ extracted from the
 semi-analytical computation (Fig.~\ref{fig} (c)). We find a good agreement, the position
maximum is slightly smaller than the position of the inflexion. Eventually, we compute
$\max_{\beta}\tau$ as a function of $\Lambda$ (Fig.~\ref{fig} (d)). We find a linear growth, in agreement with
the observation of figure~\ref{fig} (b).

Theses results all give the same information (Fig.~\ref{fig} (c,d), Eq.~(\ref{log})), they
show that to obtain a good prediction up to inverse temperature $\beta$, one needs to increase the
cut-off $\Lambda$ exponentially with $\beta$. Transposed to AMS and based on results of
figure~\ref{trans_AMS} (d), this means that when using the
linear reaction coordinate one needs to increase $N$ exponentially with $\beta$ in order to
have a given precision on the estimate of durations. An exponential behavior of $N$ translates into a
logarithmic convergence. This is consistent with the numerical convergence study.

\section{Discussion and conclusion \label{conc}}
This numerical work has focused on a detailed statistical behavior of AMS on simple 1-D and 2-D examples.
We have emphasized the difference between known theoretical results and the statistics observed in
our models in dimension larger than 1. Up to now, theoretical results have been demonstrated only in
dimension 1 or when a perfect version of the algorithm is achieved. The discrepancies observed in 1-D models
are in fact related to the discretised dynamics itself. This framework has not been
considered theoretically either.
There is therefore no contradictions with our numerical results.
The key result of this paper is the sensitivity of the statistics to the choice of a reactive coordinate
in dimension larger than 1.

In particular, the occurrence of a phase transition leads to nonstandard and unexpected statistics
for either the number of iterations (and thus \emph{a posteriori} estimates of the crossing probability)
or the duration of reactive trajectories.
This is revealed by the triple-well model during phase transition. In this case, the
choice of a poor reaction coordinate, namely a linear coordinate here, yields
deviations from pure Poissonian law in the asymptotic regime $N \to \infty$.
It is very likely that this type of critical behaviour will be found in more complex and realistic systems.
The algorithm can be very helpful for detecting an improper reaction coordinate, which will
provide transition paths which are not the most probable ones.

The approach taken for the three-level Markov model could in principle be applied to
more complex systems, with more intermediate states. The logarithmic convergence for the duration of reactive
trajectories is likely to be a generic situation as well.
The analysis and the modeling proposed can help shed light on
the quality of the results when the convergence of the algorithm is unsure.

All things considered, we found that the implementation and test of the algorithm stressed the
different convergence regimes which can be observed.
In particular, it calls for a careful examination of the types of asymptotic $N\rightarrow \infty$
behaviour of the algorithm. Our analysis shows that it depends weakly on the model and
strongly on the reaction coordinate. This is necessary to provide intervals of confidence valid
in all cases, in particular in those where a general convergence study is out of the question, but
the properties of available reaction coordinate can be easily tested. Note that reaction coordinate
selecting improbable pathways are not altogether useless: they allow one to sample a
different kind of physics. However, one must be wary that in that case, AMS does not compute the
typical behaviour of the model.

Nevertheless, the algorithm appears to be a promising tool for the study of rare events in spatially extended systems, such as those seen in the 2D Navier--Stokes equation~\cite{BS}. The first steps toward the understanding of these system can consist in the study of a typical metastable SPDE, such as the Ginzburg--Landau equation \cite{GZ}. Meanwhile, the algorithm also proves itself efficient at investigating non-equilibrium model such as the AB model \cite{bt_niso,ABnum}

\section*{acknowledgments}
The authors thank an anonymous reviewer for helpful comments that helped improve the manuscript. J. Rolland acknowledges the hospitality of the IAU of Frankfurt Goethe University where the modifications where implemented.

\end{document}